\documentclass{article}[12pt]
\usepackage[a4paper, total={6in, 9in}]{geometry}

\usepackage{amsmath,amssymb,dsfont,amsfonts,amsthm,amscd}
\usepackage{tabularx}
\usepackage{graphicx}
\usepackage{pictexwd} 
\usepackage{times}
\usepackage{fullpage}
\usepackage{xcolor}
\usepackage{color}
\usepackage{float}
\usepackage[pdfencoding=auto]{hyperref}
\usepackage{url}
\usepackage{array}
\usepackage{bm}
\usepackage{authblk}
\usepackage{natbib}
\usepackage{setspace}
\usepackage{booktabs}
\usepackage{stackengine}
\usepackage{rotating}
\usepackage[caption=false]{subfig}
\stackMath

\DeclareMathOperator*{\argmax}{arg\,max}
\DeclareMathOperator*{\argmin}{arg\,min}

\title{Counterfactual Analysis and Target Setting in Benchmarking}

\author[1]{Peter Bogetoft\thanks{Peter Bogetoft: \tt{pb.eco@cbs.dk}}}
\author[2]{Jasone Ram\'{\i}rez-Ayerbe\thanks{Jasone Ram\'{\i}rez-Ayerbe: \tt{mrayerbe@us.es}}}
\author[1]{Dolores Romero Morales\thanks{Dolores Romero Morales: \tt{drm.eco@cbs.dk}}}
\affil[1]{Department of Economics, Copenhagen Business School, Frederiksberg, Denmark}
\affil[2]{Instituto de Matem\'aticas de la Universidad de Sevilla, Seville, Spain}
\date{}

\providecommand{\keywords}[1]{\textbf{\textit{Keywords---}} #1}

\usepackage{pdfpages}

\begin{document}

\includepdf{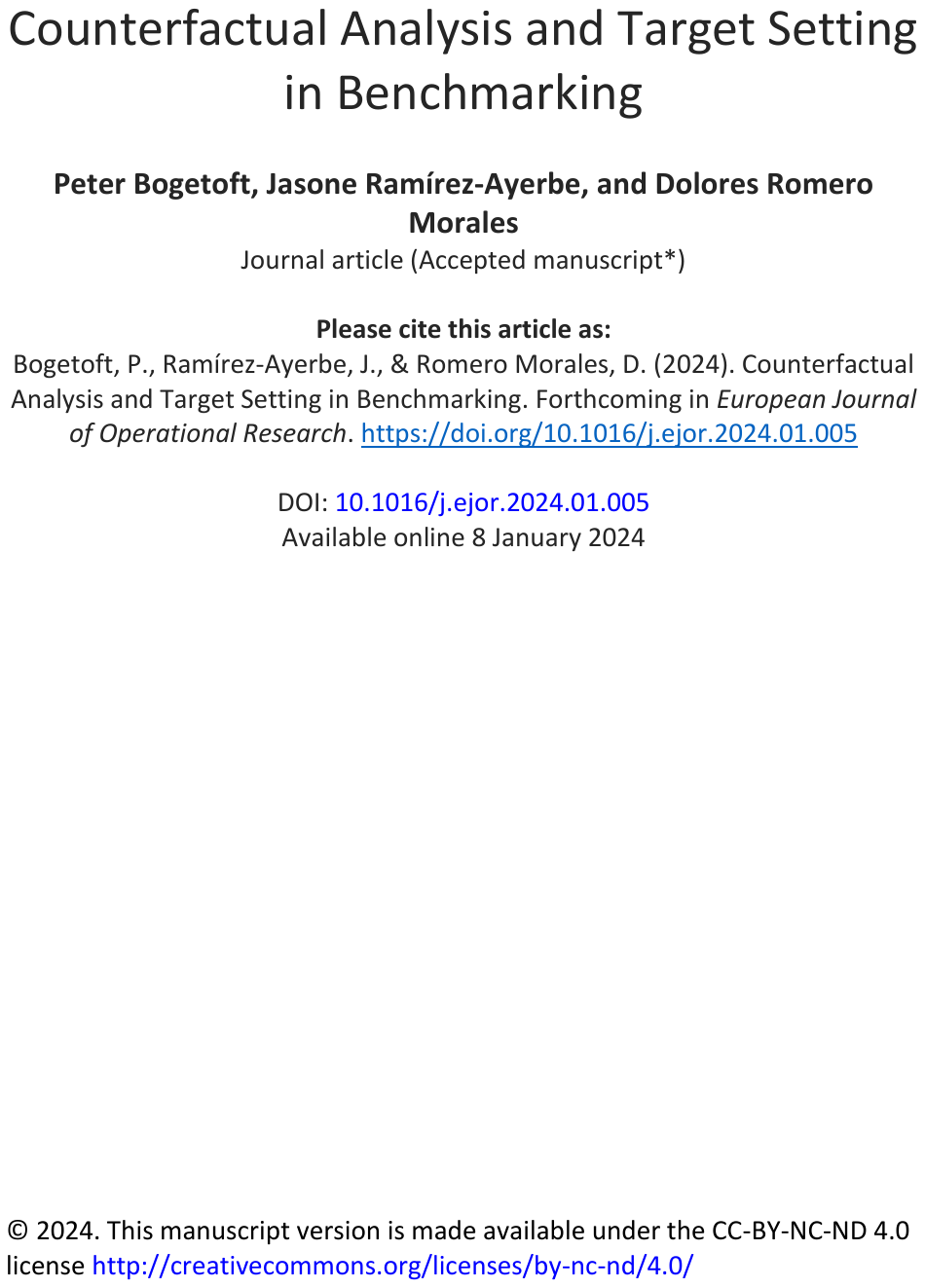}

\newpage 

\maketitle

\begin{abstract}
	
	Data Envelopment Analysis (DEA) allows us to capture the complex relationship between multiple inputs and outputs in firms and organizations. Unfortunately, managers may find it hard to understand a DEA model and this may lead to mistrust in the analyses and to difficulties in deriving actionable information from the model. In this paper, we propose to use the ideas of target setting in DEA and of counterfactual analysis in Machine Learning to overcome these problems. We define DEA counterfactuals or targets as alternative combinations of inputs and outputs that are close to the original inputs and outputs of the firm and lead to desired improvements in its performance. We formulate the problem of finding counterfactuals as a bilevel optimization model. For a rich class of cost functions, reflecting the effort an inefficient firm will need to spend to change to its counterfactual, finding counterfactual explanations boils down to solving Mixed Integer Convex Quadratic Problems with linear constraints. We illustrate our approach using both a small numerical example and a real-world dataset on banking branches.
	
\end{abstract}

\keywords{Data Envelopment Analysis; Benchmarking; DEA Targets; Counterfactual Explanations; Bilevel Optimization}

\section{Introduction}
\label{sec:intro}

In surveys among business managers, benchmarking is consistently ranked as one of the most popular management tools \citep{ManagemenTools2015,ManagementToolsTrends}. The core of benchmarking is relative performance evaluation. The performance of one entity is compared to that of a group of other entities. The evaluated ``entity'' can be a firm, organization, manager, product or process. In the following, it will be referred to simply as a Decision Making Unit (DMU).

There are many benchmarking approaches and they can serve different  purposes, such as, facilitating \emph{learning}, \emph{decision making} and \emph{incentive design}. Some approaches are very simple and rely on the comparison of a DMU's Key Performance Indicators (KPIs) to those of a selected peer group of DMUs. These KPIs are basically partial productivity measures (e.g., labour productivity, yield per hectare, etc.). This makes KPI based benchmarking easy to understand, but also potentially misleading by ignoring the role of other inputs and outputs in real DMUs. More advanced benchmarking approaches rely on frontier models using mathematical programming, e.g., Data Envelopment Analysis (DEA), and Econometrics, e.g., Stochastic Frontier Analysis (SFA), and they allow us to explicitly model the complex interaction between the multiple inputs and outputs among best-practice DMUs, cf.\ e.g.~\cite{BogetoftOtto2011,Coopereabook1995,parmeter2019combining,zhu_book2016}.

In this paper we focus on DEA based benchmarking. To construct the best practice performance frontier and evaluate the efficiency of a DMU relative to this frontier, DEA introduces a minimum of production economic regularities, typically convexity, and uses linear or mixed integer programming to capture the relationship between multiple inputs and outputs of a DMU. In this sense, and in the eyes of the modeller, the method is well-defined and several of the properties of the model will be understandable from the production economic regularities. Still, from the point of view of the evaluated DMUs, the model will appear very much like a black box. Understanding a multiple input and multiple output structure is basically difficult. Also, in DEA, there is no explicit formula showing the impact of specific inputs on specific outputs as in SFA or other econometrics based approaches. This has led some researchers to look for extra information and structure of DEA models, most notably by viewing the black box as a network of more specific process, cf.\ e.g.~\cite{Cherchye_ea2013,FareGrosskopf2000,Kao2009}.

The black box nature of DEA models may lead to some algorithm aversion and mistrust in the model, and to difficulties in deriving actionable information from the model beyond the efficiency scores. To overcome this and to get insights into the functioning of a DEA model, there are several strands of literature and tools that can be useful. The Multiple Criteria Decision Making (MCDM) literature has developed several ways in which complicated sets of alternatives can be explored and presented to a decision maker.  Also, in DEA, there is already a considerable literature on finding targets that a firm can investigate in attempts to find attractive alternative production plans. Last, but not least, it may be interesting to look for counterfactual explanations much like they are used in machine learning.

In this paper, we propose the use of counterfactual and target analyses to understand and explain the efficiencies of individual DMUs, to learn about the estimated best practice technology, and to help answer what-if questions that are relevant in operational, tactical and strategic planning efforts \citep{Bogetoft2012}. In a DEA context, counterfactual and target analyses can help with learning, decision making and incentive design. In terms of learning, the DMU may be interested to know what simple changes in features (inputs and outputs) lead to a higher efficiency level. In the application investigated in our numerical section, this can be, for instance, how many credit officers or tellers a bank branch should remove to become fully efficient. This may help the evaluated DMU learn about and gain trust in the underlying modelling. In terms of decision making, targets and counterfactual explanations may help guide the decision process by offering the smallest, the most plausible and actionable, and the least costly changes that lead to a desired boost in performance. It depends on the context how to define the least costly, or the most plausible or actionable improvement paths. In some cases it may be easier to reduce all inputs more or less the same (lawn mowing), while in other cases certain inputs should be reduced more aggressively than others, cf.~\cite{AntleBogetoft2019}. Referring back to the application in the numerical section, reducing the use of different labor inputs could for example take into account the power of different labor unions and the capacity of management to struggle with multiple employee groups simultaneously.  Lastly, targets and counterfactual explanations may be useful in connection with incentive provisions. DEA models are routinely used by regulators of natural monopoly networks to incentivize cost reductions and service improvement, cf.\ e.g. \cite{Haney2009} and later updates in \cite{AgrellBogetoftHANDBOOK, Bogetoft2012}. Regulated firms will naturally look for the easiest way to accommodate the regulator's efficiency thresholds. Counterfactual explanations may in such cases serve to guide the optimal strategic responses to the regulator's requirements.

Unfortunately, it is not an entirely trivial task to properly determine targets and construct counterfactual explanations in a DEA context. We need to find alternative solutions that are in some sense close to the existing input-output combination used by a DMU. This involves finding ``close'' alternatives in the complement of a convex set \citep{thach1988design}. In this paper, we investigate different ways to measure the closeness between a DMU and its counterfactual DMU, or the cost of moving from an existing input-output profile to an alternative target. In particular, we suggest to use combinations of  $\ell_0$, $\ell_1$, and  $\ell_2$ norms. We also consider both changes in input and output features and show how to formulate the problems in DEA models with different returns to scale assumptions. We show how determining targets and constructing counterfactual explanations leads to a bilevel optimization model, that can be reformulated as a Mixed Integer Convex Quadratic Problem with linear constraints. We illustrate our approach on both a small numerical example as well as a large scale real-world dataset involving bank branches.

The outline of the paper is as follows. In Section \ref{sec:litrev} we review the relevant literature. In Section \ref{sec:setting} we introduce the necessary DEA notation for constructing targets and counterfactual explanations, as well as a small numerical example. In Section \ref{sec:bilevel} we describe our bilevel optimization formulation and its reformulation as a Mixed Integer Convex Quadratic Problem with linear constraints.  In Section \ref{sec:num} we illustrate our approach with real-world data on bank branches. We end the paper with conclusions in Section \ref{sec:conclusion}. In the Appendix, we extend the analysis by investigating alternative returns to scale and by investigating changes in the outputs rather than the inputs.

\section{Background and Literature}
\label{sec:litrev}

In this section, we give some background on DEA benchmarking, in particular on directional and interactive benchmarking, on target setting in DEA, and on counterfactual analysis from interpretable machine learning.

Data Envelopment Analysis, DEA, was first introduced in \cite{Charnes1978, Charnes1979} as a tool for measuring efficiency and productivity of decision making units, DMUs. The idea of DEA is to model the production possibilities of the DMUs and to measure the performance of the individual DMUs relative to the production possibility frontier. The modelling is based on observed practices that form activities in a Linear Programming (LP) based activity analysis model. 

Most studies use DEA models primarily to measure the relative efficiency of the DMUs. The benchmarking framework, the LP based activity analysis model, does however allow us to explore a series of other questions.  In fact, the benchmarking framework can serve as a learning lab and decision support tool for managers. In the DEA literature, this perspective has been emphasized by the idea of  interactive benchmarking. Interactive benchmarking and associated easy to use software has been used in a series of applications and consultancy projects, cf.\ e.g.~\cite{Bogetoft2012}. The idea is that a DMU can search for alternative and attractive production possibilities and hereby learn about the technology, explore possible changes and trade-offs and look for least cost changes that allow for necessary performance improvements, cf.\ also our discussion of learning, decision making and incentive and regulation applications in the introduction.

One way to illustrate the idea of interactive benchmarking is as in Figure \ref{fig:Directional} below. A DMU has used two inputs to produce two outputs. Based on the data from other DMUs, an estimate of the best practice technology has been established as illustrated by the piecewise linear input and output isoquants. The DMU may now be interested in exploring alternative paths towards best practices.  One possibility is to save a lot of input 2 and somewhat less of input 1, i.e., to move in the direction $\bm{d}_x$ illustrated by the arrow in the left panel. If the aim is to become fully efficient, this approach suggests that the DMU instead of the present (input,output) combination $(\bm{x},\bm{y})$ should consider the alternative $(\bm{\hat{x}},\bm{y})$. A similar logic could be used on the output side keeping the inputs fixed as illustrated in the right panel where we assume that more of a proportional increase in the two outputs is strived at. Of course, in reality, one can combine also changes in the inputs and outputs.

\begin{figure}[htb]
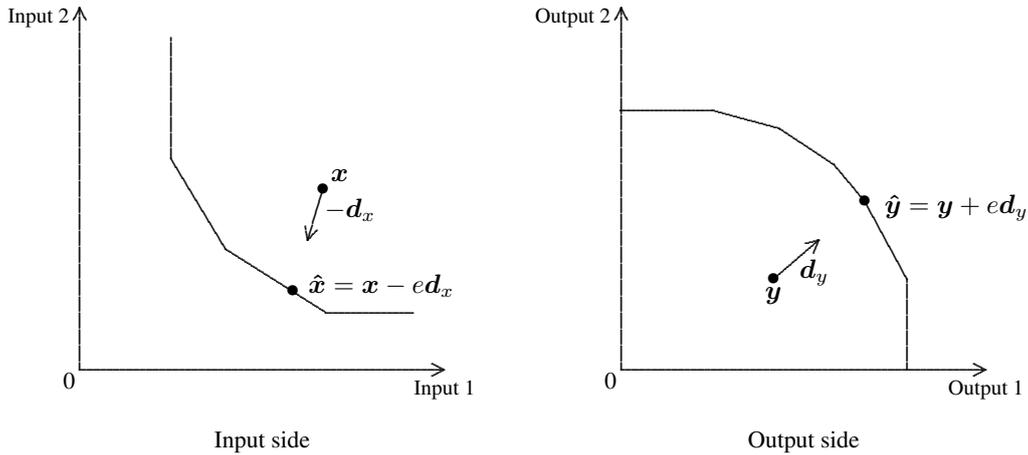

	\begin{center}
		\begin{minipage}{.45\linewidth}
			\beginpicture
			\setcoordinatesystem units <.8mm, .8mm>
			\normalgraphs
			\arrow <5pt> [.5,1] from 0 0 to 0 60
			\arrow <5pt> [.5,1] from 0 0 to 60 0
			\put {\footnotesize Input 1}  [tc] <0pt,-4pt> at 60 0
			\put {\footnotesize Input 2}  [tr] <-4pt,0pt> at 0 60
			\put {\small 0}  [cc] <-4pt,-4pt> at 0 0
			\linethickness=5pt
			\setlinear
			\plot 15 55  15 35  24 20  40.5 9.5   55 9.5  /
			\put {$\bullet$} at 40 30
			\put {$\bullet$} at 35 13
			\put {$\bm{x}$} [Bl] <3pt,3pt>  at 40 30
			\put {$\bm{\hat{x}}=\bm{x}-e\bm{d}_x$} [bl] <6pt,-2pt>  at 35 13
			\arrow <5pt> [.5,1] from 40 30 to 37.5 21.5
			\put {$-\bm{d}_x$} [Bl] <3pt,3pt>  at 39 24
			\put {\small Input side} [c] at 30 -12
			\endpicture
		\end{minipage}
		\qquad
		\begin{minipage}{.5\linewidth}
			\beginpicture
			\setcoordinatesystem units <.8mm, .8mm>
			\normalgraphs
			\arrow <5pt> [.5,1] from 0 0 to 0 60
			\arrow <5pt> [.5,1] from 0 0 to 60 0
			\put {\footnotesize Output 1}  [tc] <0pt,-4pt> at 60 0
			\put {\footnotesize Output 2}  [tr] <-4pt,0pt> at 0 60
			\put {\small 0}  [cc] <-4pt,-4pt> at 0 0
			\linethickness=5pt
			\setlinear
			\plot 47 0   47 15   40 28   35 34  26 40  15 43  0 43  /
			\put {$\bullet$} at  25 15
			\put {$\bullet$} at  40 28
			\put {$\bm{y}$} [tc] <0pt,-3pt>  at 25 15
			\put {$\bm{\hat{y}}=\bm{y}+e\bm{d}_y$} [Bl] <3pt,1pt>  at 42 25.5
			\arrow <5pt> [.5,1] from 25 15 to 32.5 21.5
			\put {$\bm{d}_y$} [Bl] <3pt,3pt>  at 28 14
			\put {\small Output side} [c] at 30 -12
			\endpicture
		\end{minipage}
	\end{center}
	\vskip1ex
	\caption{Directional search for an alternative production plan to $(\bm{x},\bm{y})$ along $(\bm{d}_x,\bm{d}_y)$ using (DIR) }\label{fig:Directional}
\end{figure}

Formally, the directional distance function approach, sometimes referred to as the excess problem, requires solving the following mathematical programming problem
\begin{align}
	\label{eq:directional}\tag{DIR}
	\max \quad & \{e \: \mid \:  (\bm{x}-e\bm{d}_x, \bm{y}+e\bm{d}_y) \in T^*\},
\end{align}
where $\bm{x}$ and $\bm{y}$ are the present values of the inputs and output vectors, $\bm{d}_x$ and $\bm{d}_y$ are the improvement directions in input and output space, $T^*$ is the estimated set of feasible (input,output) combinations, and $e$ is the magnitude of the movement.

In the DEA literature, the direction $(\bm{d}_x,\bm{d}_y)$ is often thought as parameters that are given and the excess as one of many possible ways to measure distance to the frontier. A few authors have advocated that some directions are more natural than others and there have been attempts to endogenize the choice of this direction, cf.\ e.g.~\cite{Bogetoft1999, FareGrosskopfWittaker2013, Fare_ea2017,petersen2018directional, ZofioPastorAparicio2013}. One can also think of the improvement directions as reflecting the underlying strategy of the DMU or simply as a steering tool that the DMU uses to create one or more interesting points on the frontier.

\begin{figure}[h]
	\begin{center}
		\includegraphics[scale=0.6]{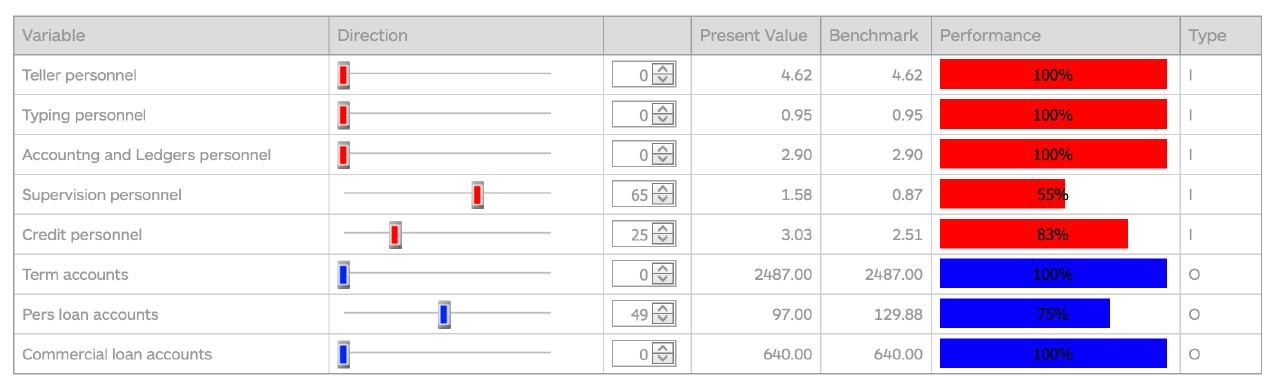}
		\caption{Directional search in Interactive Benchmarking software. Real-world dataset of bank branches in Section \ref{sec:num}}
		\label{fig:IB_example}
	\end{center}
\end{figure}

Figure \ref{fig:IB_example} illustrates the real-world example involving bank branches from the application section. The analysis is here done using the directional distance function approach (DIR) as implemented in the so-called Interactive Benchmarking software, cf.\ \cite{Bogetoft2012}. The search ``Direction" is chosen by adjusting the horizontal handles for each input and output and is expressed in percentages of the existing inputs and outputs. The resulting best practice alternative is illustrated in the ``Benchmark'' column. We see that the DMU in this example expresses an interest in reducing Supervision and Credit personnel but simultaneously seeks to increase the number of personal loan accounts.

Applications of interactive benchmarking have typically been in settings where the DMU in a trial-and-error like process seeks alternative production plans. Such processes can certainly be useful in attempts to learn about and gain trust in the modelling, to guide decision making and to find the perfomance enhancing changes that a DMU may find relatively easy to implement.  From the point of view of Multiple Criteria Decision Making (MCDM) we can think of such processes as based on progressive articulation of preferences and alternatives, cf.\ e.g.\ the taxonomy of MCDM methods suggested in \cite{Rostamy_ea2017}.

It is clear from this small example, however, that the use of an interactive process guided solely by the DMU may not always be the best approach. If there are more than a few inputs and outputs, the process can become difficult to steer towards some underlying optimal compromise between the many possible changes in inputs and outputs. In such cases, the so-called prior articulation of preferences may be more useful. If the DMU can express its preferences for different changes, e.g., as a cost of change function $C((\bm{x},\bm{y}),(\bm{x}^*,\bm{y}^*))$ giving the cost of moving from the present production plan $(\bm{x},\bm{y})$ to any new production plan $(\bm{x}^*,\bm{y}^*)$, then a systematic search for the optimal change is possible. The approach of this paper is based on this idea. We consider a class of cost functions and show how to find optimal changes in inputs and outputs using bilevel optimization. In this sense, it corresponds to endogenizing the directional choice so as to make the necessary changes in inputs and outputs as small as possible. Of course, by varying the parameters of the cost function, one can also generate a reasonably representative set of alternative production plans that the DMU can then choose from. This would correspond to the idea of a prior articulation of alternatives approach in the MCDM taxonomy.

The idea of introducing a cost of change function to guide the search for alternative production plans is closely related to the idea of targets in DEA. At a general level, a target is here understood as an alternative production plan that a DMU should move to.\footnote{In \cite{AparicioRuizSirvent2007}, the authors distinguish between setting targets and benchmarking in the sense that targets are the coordinates of a projection point, which is not necessarily an observed DMU, whereas benchmarks are real observed DMUs. This distinction can certainly be relevant is several contexts, but it is not how we use targets here. We use benchmarking as the general term for relative performance comparison, and target to designate an alternative production plan that a DMU should choose to improve performance in the easiest possible way.} There has been a series of interesting DEA papers on the determination of targets using the principle of least action, see for example \cite{ Aparicio_et_al.2014c} and the references in here. In \cite{ Aparicio_et_al.2014c}, the authors explicitly introduce  the principle of least action referring to the idea in physics that nature always finds the most efficient course of action. The general argument underlying these approaches is that an inefficient firm should achieve technical efficiency with a minimum amount of effort. Different solutions have been proposed using different distance measures or what we call the cost of change. In many papers, this corresponds to minimizing the distance to the efficient frontier in contrast to the traditional efficiency measurement problem, where we are looking for the largest possible savings or the largest possible expansions of the services provided. A good example is \cite{ AparicioPastor2013}. In this sense, our idea of finding close counterfactuals fits nicely into the DEA literature.

The choice of targets has also been discussed in connection with the slack problem in radial DEA measures. A Farrell projection may not lead to a Pareto efficient point and in a second stage, it is therefore common to discuss close alternatives that are fully Pareto efficient. Again, different solutions -- using for example constraints or the determination of all full facets -- have been proposed, cf.~\cite{ Aparicio_et_all_2017Omega}. Identifying all facets and measuring distance to these like in \cite{Portela2003} is theoretically attractive but computationally cumbersome in most applications.

Our approach can be seen as a generalization of the literature on targets. In particular, if $E^*=1$ and the considered cost function coincides with a norm or with a typical technical efficiency measure, then the previous DEA target approaches are particular cases of the general approach introduced in this paper. We formulate the target setting problem for general cost-of-change functions using a bilevel program\footnote{The idea of using a bilevel linear programming approach has also appeared in the DEA literature. It should be noted in particular that \cite{ Aparicio_et_all_2017Omega} proposed to resort to a bilevel linear programming model when strictly efficient targets are to be identified using the Russel output measure to capture cost-of-change.}, and we reformulate the constraints to get tractable mathematical optimization problems. Using combinations of $\ell_0$,  $\ell_1$ and $\ell_2$ norms, as we do in the illustrations, the resulting problems are Mixed Integer Convex Quadratic Problems with linear constraints. It is worthwhile to note also, that we do not necessarily require the target to be Pareto efficient, allowing for the possibility that a DMU may not seek to become fully efficient but, for example, just 90\% Farrell efficient which also implies that targets may be on non-full facets.

In interpretable machine learning \citep{du2019techniques,rudin2022interpretable}, counterfactual analysis is used to explain the predictions made for individual instances \citep{guidotti2022counterfactual,karimi2022survey,wachter2017counterfactual}.  Machine learning approaches like Deep Learning, Random Forests, Support Vector Machines, and XGBoost are often seen as powerful tools in terms of learning accuracy but also as black boxes in terms of how the model arrives at its outcome. Therefore, regulations from, among others the EU, are enforcing more transparency in the so-called field of algorithmic decision making \citep{EUwhitepaperAI20,goodmanAIM17}. There is a paramount of tools being developed in the nascent field of explainable artificial intelligence to help understand how tools in machine learning and artificial intelligence make decisions \citep{lundberg2017unified,martensMISQ14,molnar2020interpretable}. The focus of this paper is on counterfactual analysis tools. The starting point is an individual instance for which the model predicts an undesired outcome. In counterfactual analysis, one is interested in building an alternative instance, the so-called counterfactual instance, revealing how to change the features of the current instance so that the model predicts a desired outcome for the counterfactual instance. The counterfactual explanation problem is written as a mathematical optimization problem. To define the problem, one needs to model the feasible space, a cost function measuring the cost of the movement from the current instance to the counterfactual one, and a set of constraints that ensures that the counterfactual explanation is predicted with the desired outcome. In general, the counterfactual explanation problem reads as a constrained nonlinear problem but, for score-based classifiers and cost functions defined by a convex combination of the norms $\ell_0$,  $\ell_1$ and $\ell_2$, equivalent Mixed Integer Linear Programming or Mixed Integer Convex Quadratic with Linear Constraints formulations can be defined, see, e.g., \cite{carrizosaESWA24,Fischetti2018,parmentier2021optimal}.

In the following, we combine the ideas of DEA, least action targets, and counterfactual explanations. We formulate and solve bilevel optimization models to determine ``close'' alternative production plans or counterfactual explanations in DEA models that lead to desired relative performance levels and also take into account the strategic preferences of the entity.

\section{The Setting}
\label{sec:setting}

We consider $K+1$ DMUs (indexed by $k$), using $I$ inputs, $\bm{x}^k=(x^k_1,\dots,x_I^k)^\top \in \mathbb{R}_{+}^I$, to produce $O$ outputs, $\bm{y}^k=(y_1^k,\dots,y_O^k)^\top \in \mathbb{R}_{+}^O$. Hereafter, we will write $(\bm{x}^k, \bm{y}^k)$ to refer to production plan of DMU $k$, $k=0,1,\dots,K$.

Let $T$ be the technology set, with
$$
T=\{(\bm{x},\bm{y}) \in \mathbb{R}_{+}^I \times \mathbb{R}_+^O \quad | \quad \bm{x} \text{ can produce } \bm{y} \}.
$$
We will initially estimate $T$ by the classical DEA model. It determines the empirical reference technology $T^{*}$ as the smallest subset of  $\mathbb{R}_{+}^I \times \mathbb{R}_+^O$ that contains the actual $K+1$ observations, and satisfies the classical DEA regularities of convexity, free-disposability in inputs and outputs, and Constant Returns to Scale (CRS). It is easy to see that the estimated technology can be described as:
$$
T^{*}(\text{CRS})=\{ (\bm{x},\bm{y}) \in \mathbb{R}_{+}^I \times \mathbb{R}_+^O \quad | \quad \exists \bm{\lambda} \in \mathbb{R}_+^{K+1}: \bm{x} \geq \sum_{k=0}^K \lambda^k \bm{x}^k, \bm{y}\leq \sum_{k=0}^K \lambda^k \bm{y}^k \}.
$$

To measure the efficiency of a firm, we will initially use the so-called Farrell input-oriented efficiency. It measures the efficiency of a DMU, say DMU $0$, as the largest proportional reduction $E^0$ of all its inputs $\bm{x}^0$ that allows the production of its present outputs $\bm{y}^0$ in the technology $T^{*}$. Hence, it is equal to the optimal solution value of the following LP formulation
\begin{align}
	\label{eq:primal}\tag{DEA}
	\min_{E,\lambda^0,\dots,\lambda^K} \quad & E \\
	\text{s.t.} \quad & E \bm{x}^0 \geq \sum_{k=0}^K \lambda^k \bm{x}^k \nonumber\\
	& \bm{y}^0 \leq \sum_{k=0}^K \lambda^k \bm{y}^k \nonumber\\
	& 0 \leq E \leq 1 \nonumber\\
	& \bm{\lambda} \in \mathbb{R}^{K+1}_{+}. \nonumber
\end{align}
This DEA model has $K+2$ decision variables, $I$ linear input constraints and $O$ linear output constraints. Hereafter, we will refer to the optimal objective value of (DEA), say $E^0$, as the efficiency of DMU 0.

In the following, and assuming that firm 0 with production plan $(\bm{x}^0,\bm{y}^0)$ is not fully efficient, $E^0<1$, we will show how to calculate a \emph{counterfactual explanation} with a desired efficiency level $E^* > E^0$, i.e., the minimum changes needed in the inputs of the firm, $\bm{x}^0$, in order to obtain an efficiency $E^{*}$. Given a cost function $C(\bm{x}^0,\hat{\bm x})$ that measures the cost of moving from the present inputs $\bm{x}^0$ to the new counterfactual inputs $\hat{\bm x}$, and a set $\mathcal{X}(\bm{x}^0)$ defining the feasible space for $\hat{\bm x}$, the counterfactual explanation for $\bm{x}^0$ is found solving the following optimization problem:
\begin{align*}
	\min_{\hat{\bm x}} \quad &C(\bm{x}^0,\hat{\bm x})\\
	\text{s.t.} \quad & \hat{\bm x} \in \mathcal{X}(\bm{x}^0)\\
	& (\hat{\bm x}, \bm{y}^0) \quad \text{has at least an efficiency of $E^{*}$}.
\end{align*}
With respect to $C(\bm{x}^0,\hat{\bm x})$, different norms can be used to measure the difficulty of changing the inputs. A DMU may, for example, be interested to minimize the sum of the squared deviations between the present and the counterfactual inputs. We model this using the squared Euclidean norm $\ell_2^2$. Likewise, there may be an interest in minimizing the absolute value of the deviations, which we can proxy using the  $\ell_1$ norm, or the number of inputs changed, which we can capture with the  $\ell_0$ norm. When it comes to $\mathcal{X}(\bm{x}^0)$, this would include the nonnegativity of $\hat{\bm x}$, as well as domain knowledge specific constraints. With this approach, we detect the most important inputs in terms of the impact they have on the DMU's efficiency, and with enough flexibility to consider different costs of changing depending on the DMU's characteristics.

In the next section, we will show that finding counterfactual explanations involves solving a bilevel optimization problem of minimizing the changes in inputs and solving the above DEA problem at the same time. In the Appendix, we will also discuss how the counterfactual analysis approach can be extended to other technologies and to other efficiency measures like the output-oriented Farrell efficiency and other DEA technologies.

Before turning to the details of the bilevel optimization problem, it is useful to illustrate the idea of counterfactual explanations using a small numerical example. Suppose we have four firms with the inputs, outputs, and Farrell input efficiencies as in Table \ref{tab:example}. The efficiency has been calculated solving the classical DEA model with CRS, namely \eqref{eq:primal}. In this example, firms 1 and 2 are fully efficient, whereas firms 3 and 4 are not.
\begin{table}[h]
	\centering
	\begin{tabular}{l|cccc}
		Firm & $x_1$ & $x_2$ & $y$ & $E$\\ \hline
		1 & 0.50 & 1 & 1& 1\\
		2&1.50 & 0.50 & 1 & 1\\
		3&1.75& 1.25& 1& 0.59\\
		4&2.50&1.25&1&0.50\\
	\end{tabular}
	\caption{Inputs, outputs and corresponding Farrell input-efficiency of 4 different firms}
	\label{tab:example}
\end{table}

\begin{figure}[htb]
	\subfloat[Explanations for firm 3\label{fig:examplesol}]
	{\includegraphics[width=.43\linewidth]{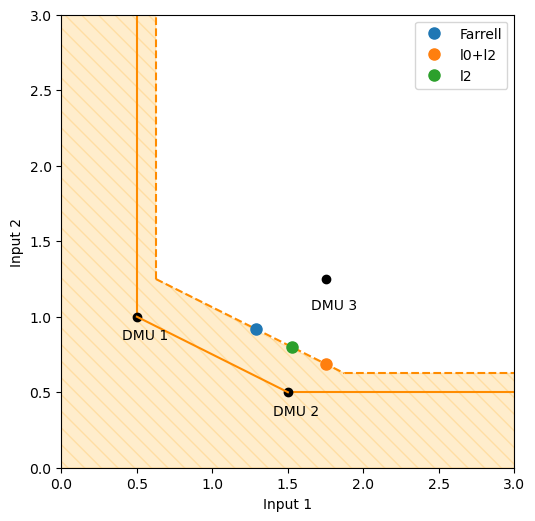}}\hfill
	\subfloat[Explanations for firm 4\label{fig:examplesol2}]
	{\includegraphics[width=.43\linewidth]{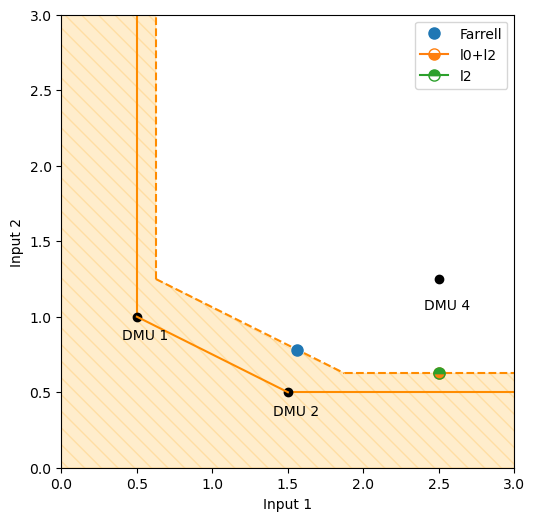}}
	\caption{Counterfactual explanations for firms 3 and 4 in Tables \ref{tab:sol} and \ref{tab:sol2} respectively imposing $E^{*}=0.8$ and different cost functions}
	\label{fig:examplesol1y2}
\end{figure}

First, we want to know the changes needed in $\bm{x}^3$ for firm 3 to have a new efficiency $E^{*}$ of at least 80\%. Since we only have two inputs, we can illustrate this graphically as in Figure \ref{fig:examplesol}. The results are shown in Table \ref{tab:sol} for different cost functions. It can be seen that we in all cases get exactly 80\% efficiency with the new inputs.  We see from column $\ell_2^2$ that the Farrell solution is further away from the original inputs than the counterfactual solution based on the Euclidean norm. To the extent that difficulties of change is captured by the $\ell_2^2$ norm, we can conclude that the Farrell solution is not ideal. Moreover, in the Farrell solution one must by definition change both inputs, see column $\ell_0$. Using a cost function combining the $\ell_0$ norm and the squared Euclidean norm, denoted by $\ell_0 + \ell_2$, one penalizes the number of inputs changed. With this we detect the one input that should change in order to obtain a higher efficiency, namely the second input. In contexts like negotiations with various input suppliers, it is often more practical to focus negotiations on just one or a select few inputs, instead of dealing with all inputs at the same time.
\begin{table}[h]
	\centering
	\begin{tabular}{l|cccccc}
		Cost function & $\hat{x}_1$ & $\hat{x}_2$ & $y$ & $E$ & $\ell_2^2$ & $\ell_0$\\ \hline
		Farrell &1.29& 0.92 & 1& 0.8&0.32&2\\
		$\ell_0 + \ell_2$&1.75 &0.69 & 1 & 0.8&0.31&1\\
		$\ell_2$&1.53& 0.80& 1& 0.8&0.25&2
	\end{tabular}
	\caption{Counterfactual explanations for firm 3 in Table \ref{tab:example} imposing $E^{*}=0.8$ and different cost functions}
	\label{tab:sol}
\end{table}

Let us now focus on firm 4 and again find a counterfactual instance with at least 80\% efficiency. The results are shown in Table \ref{tab:sol2} and Figure \ref{fig:examplesol2}. Notice how in the Farrell case one obtains again the farthest solution and also the least sparse from the three of them. As for the counterfactual explanations with our methodology, the inputs nearest to the original DMU that give us the desired efficiency are in a non full-facet of the efficiency frontier. By using the Farrell to measure the desired efficiency level, we only need to change one input, namely, the second input, and can have ``slack" in the first input. We here deviate from \cite{Aparicio_et_all_2017Omega}, in which the authors look for targets on the strongly efficient frontier, i.e., without slack.

\begin{table}[h]
	\centering
	\begin{tabular}{l|cccccc}
		Cost function & $\hat{x}_1$ & $\hat{x}_2$ & $y$ & E & $\ell_2^2$  &$\ell_0$\\ \hline
		Farrell & 1.56& 0.78 & 1& 0.8&1.10&2\\
		$\ell_0+\ell_2$&2.50 & 0.63 & 1 & 0.8&0.39&1\\
		$\ell_2$&2.50& 0.63& 1& 0.8&0.39&1
	\end{tabular}
	\caption{Counterfactual explanations for firm 4 in Table \ref{tab:example} imposing $E^{*}=0.8$ and different cost functions}
	\label{tab:sol2}
\end{table}

In Figure \ref{fig:examplesol1y2}, the space where we search for the counterfactual explanation is shaded. Although in these illustrations the frontier is explicitly given, in general, the frontier points are convex combinations of the observed DMUs, and obtaining the isoquants is not easy. Therefore, when finding counterfactual explanations a bilevel optimization is needed to search for ``close'' inputs in the complement of a convex set.

\section{Bilevel optimization for counterfactual analysis in DEA}
\label{sec:bilevel}

Suppose DMU 0 is not fully efficient, i.e., the optimal objective value of Problem \eqref{eq:primal} is $E^{0}<1$. In this section, we formulate the counterfactual explanation problem in DEA, i.e., the problem that calculates the minimum cost changes in the inputs $ \bm{x}^0$ that make DMU 0 have a higher efficiency. Let $\hat{\bm x}$ be the new inputs of DMU $0$ that would make it at least $E^{*}$ efficient, with $E^{*} > E^{0}$. With this, we have defined the counterfactual instance as the one obtained changing the inputs, but in the same sense, we could define it by changing the outputs. This alternative output-based problem will be studied in the Appendix.

Since the values of the inputs are to be changed, the efficiency of the new production plan $(\hat{\bm x},\bm{y}^0)$ has to be calculated using Problem \eqref{eq:primal}. The counterfactual explanation problem in DEA reads as follows:
\begin{align}
	\footnotesize
	\min_{\bm{\hat x}} \quad & C(\bm{x}^0,\bm{\hat x}) \label{eq:upperlevel}\\
	\text{s.t.} \quad & \bm{\hat x} \in \mathbb{R}^I_{+} \label{eq:upperlevelinputs}\\
	&E \geq E^{*} \label{eq:upperleveleff}\\
	& E \in \argmin_{\bar E, \lambda^0, \dots, \lambda^K} \footnotesize \left\{ \right. \bar E \quad :
	\quad\bar  E \hat{\bm{x}} \geq \sum_{k=0}^K \lambda^k \bm{x}^k,
	\bm{y}^0 \leq \sum_{k=0}^K \lambda^k \bm{y}^k,  \bar E \geq 0, \bm{\lambda} \in \mathbb{R}^{K+1}_{+} \left. \footnotesize \right\}, \label{eq:lowerlevel2}
\end{align}
where in the upper level problem in \eqref{eq:upperlevel} we minimize the cost of changing the inputs for firm $0$, $ \bm{x}^0$, to $\hat{\bm x}$, ensuring nonnegativity of the inputs, as in constraint \eqref{eq:upperlevelinputs}, and that the efficiency is at least $E^{*}$, as in constraint \eqref{eq:upperleveleff}. The lower level problem in \eqref{eq:lowerlevel2} ensures that the efficiency of $(\hat{\bm x},\bm{y}^0)$ is correctly calculated. Therefore, as opposed to counterfactual analysis in interpretable machine learning, here we are confronted with a bilevel optimization problem. Notice also that to calculate the efficiency in the lower level problem in \eqref{eq:lowerlevel2}, the technology is already fixed, and the new DMU $(\hat{\bm x},\bm{y}^0)$ does not take part in its calculation.

In what follows, we reformulate the bilevel optimization problem \eqref{eq:upperlevel}-\eqref{eq:lowerlevel2} as a single-level model, by exploiting the optimality conditions for the lower-level problem. This can be done for convex lower-level problems that satisfy Slater's conditions, e.g., if our lower-level problem was linear. In our case, however, not all the constraints are linear, since in \eqref{eq:lowerlevel2} we have the product of decision variables $\bar E \hat{\bm{x}}$. To be able to handle this, we define new decision variables, namely, $F=\frac{1}{E}$ and $\beta^k=\frac{\lambda^k}{E}$, for $k=0,\ldots,K$. Thus, \eqref{eq:upperlevel}-\eqref{eq:lowerlevel2} is equivalent to:
\begin{align}
	\footnotesize
	\min_{\bm{\hat x},F} \quad & C(\bm{x}^0,\bm{\hat x}) \nonumber \\
	\text{s.t.} \quad & \bm{\hat x} \in \mathbb{R}^I_{+} \nonumber \\
	&F \leq F^{*} \nonumber \\
	& F \in \argmax_{\bar F, \bm{\beta}} \footnotesize \left\{ \right. \bar F \quad :
	\quad \hat{\bm{x}} \geq \sum_{k=0}^K \beta^k \bm{x}^k,
	\bar F \bm{y}^0 \leq \sum_{k=0}^K \beta^k \bm{y}^k,  \bar F \geq 0,  \bm{\beta} \in \mathbb{R}^{K+1}_{+} \left. \footnotesize \right\}.\label{eq:lowerlevelNEW2}
\end{align}

This equivalent bilevel opimization problem can now be reformulated as a single-level model. The new lower-level problem in \eqref{eq:lowerlevelNEW2} can be seen as the $\hat{\bm{x}}$-parametrized problem:
\begin{align}
	\max_{F,\bm{\beta}}\quad &  F \\
	\text{s.t.} \quad &  \hat{\bm{x}} \geq \sum_{k=0}^K \beta^k \bm{x}^k\label{eq:primm}\\
	&F \bm{y}^0 \leq \sum_{k=0}^K \beta^k \bm{y}^k \label{eq:primm2}\\
	& F \geq 0 \label{eq:primmF}\\
	& \bm{\beta} \geq \bm{0}. \label{eq:primm_l}
\end{align}

The Karush-Kuhn-Tucker (KKT) conditions, which include primal and dual feasibility, stationarity and complementarity conditions, are necessary and sufficient to characterize an optimal solution. Thus, we can replace problem \eqref{eq:lowerlevelNEW2} by its KKT conditions. Primal feasibility is given by \eqref{eq:primm}-\eqref{eq:primm_l}. Dual feasibility is given by:
\begin{gather}
	\label{eq:duall}
	\bm{\gamma}_I,\bm{\gamma}_O, \delta, \bm{\mu} \geq \bm{0},
\end{gather}
for the dual variables associated with constraints \eqref{eq:primm}-\eqref{eq:primm_l}, where $\bm{\gamma}_I \in \mathbb{R}_+^I$, $\bm{\gamma}_O \in \mathbb{R}_+^O$, $\delta \in \mathbb{R}_+$, $\bm{\mu} \in \mathbb{R}_+^{K+1}$. The stationarity conditions are as follows:
\begin{align}
	\bm{\gamma}_{O}^\top\bm{y}^{0}&-\delta = 1 \label{eq:statt}\\
	\bm{\gamma}_I^\top \bm{x}^k - \bm{\gamma}_{O}^\top \bm{y}^k& -\mu_k = 0 \quad k=0,\dots,K. \label{eq:statt2}
\end{align}
Lastly, we need the complementarity conditions for all constraints \eqref{eq:primm}-\eqref{eq:primm_l}. For constraint \eqref{eq:primm}, we have:
\begin{equation}
	\label{eq:condition1}
	\gamma_{I}^i=0 \quad \text{or} \quad \hat x_i-\sum_{k=0}^K \beta^k x_i^k=0 \quad i=1,\dots,I.
\end{equation}
In order to model this disjunction, we will introduce binary variables $u_i \in \{0,1\}$, $i=1,\dots,I$, and the following constraints using the big-M method:
\begin{equation}
	\label{eq:M1}
	\gamma_I^i \leq M_I u_i, \quad  \hat x_i-\sum_{k=0}^K \beta^k x_i^k \leq M_I(1-u_i), \quad i=1,\dots,I,
\end{equation}
where $M_I$ is a sufficiently large constant.

The same can be done for the complementarity condition for constraint \eqref{eq:primm2}, introducing binary variables $v_o\in \{0,1\}$, $o=1,\dots,O$, big-M constant $M_O$, and constraints:
\begin{equation}\label{eq:M2}
	\gamma_{O}^o \leq M_O v_o, \quad  -Fy_o^0+\sum_{k=0}^K \beta^k y_o^k \leq M_O(1-v_o), \quad o=1,\dots,O.
\end{equation}
The complementarity condition for constraint \eqref{eq:primm_l} would be the disjunction $\beta^k=0$ or $\mu_k=0$. Using the stationarity condition \eqref{eq:statt2} and again the big-M method with binary variables $w_k\in \{0,1\}$, $k=0,\dots,K$, and big-M constant $M_f$, one obtains the constraints:
\begin{equation}\label{eq:M4}
	\beta^k\leq M_f w_k, \quad \bm{\gamma}_I^\top \bm{x}^k - \bm{\gamma}_O^\top \bm{y}^k \leq M_f(1-w_k), \quad k=0,\dots,K.
\end{equation}

Finally, for constraint \eqref{eq:primmF} the complementarity condition yields $F=0$ or $\delta=0$. Remember that $F=1/E$, $0\leq E\leq 1$, thus $F$ cannot be zero by definition and we must impose $\delta=0$. Using stationarity condition \eqref{eq:statt}, this yields:
\begin{equation}
	\label{eq:dualeq}
	\bm{\gamma}_{O}^\top \bm{y}^0=1.
\end{equation}

We now reflect on the meaning of these constraints. Notice that constraints \eqref{eq:M1} and \eqref{eq:M2} model the slacks of the inputs and outputs respectively, while constraint \eqref{eq:M4} models the firms that define the frontier, i.e., the firms with which DMU 0 is to be compared. If binary variable $u_i=1$, then there is no slack in input $i$, i.e., $\hat x_i=\sum_{k=1}^K \beta^k x_i^k$, whereas if $u_i=0$ that means there is. The same happens with binary variable $v_o$, namely, it indicates whether there is a slack in output $o$. On the other hand, when $w_k=1$, then the equality of the dual constraint will hold $\bm{\gamma}_I^\top \bm{x}^k = \bm{\gamma}_{O}^\top \bm{y}^k$, i.e., firm $k$ is fully efficient and it is used to define the efficiency of the counterfactual instance. If $w_k=0$ then $\beta_k=0$, and firm $k$ is not being used to define the efficiency of the counterfactual instance. Let us go back to the example in the previous section with four firms with 2 inputs and 1 output and several choices of cost function $C$ of changing the inputs. When $C=\ell^2_2$, we can see that firm 3 is compared against firms 1 and 2, while firm 4 is compared against firm 2 only.

Notice that $\mu_k$ is only present in \eqref{eq:statt2}, thus it is free. In addition, we know that $\delta=0$. Therefore, we can transform the stationarity conditions \eqref{eq:statt} and \eqref{eq:statt2} to
\begin{align}
	\bm{\gamma}_{O}^\top\bm{y}^{0}& = 1 \label{eq:stattf}\\
	\bm{\gamma}_I^\top \bm{x}^k - \bm{\gamma}_{O}^\top \bm{y}^k& \geq 0 \quad k=1,\dots,K \label{eq:stattf2}\\
	\bm{\gamma}_I,\bm{\gamma}_O&\geq \bm{0},\label{eq:stattf3}
\end{align}
that are exactly the constraints in the dual DEA model for the Farrell output efficiency.

The new reformulation of the counterfactual explanation problem in DEA is as follows:
\begingroup
\allowdisplaybreaks
\begin{align*}
	\min_{\bm{\hat x},F, \bm{\beta},\bm{\gamma}_I,\bm{\gamma}_O,\bm{u},\bm{v},\bm{w} } \quad & C(\bm{x}^0,\bm{\hat x}) \\
	\text{s.t.} \quad 	&F\leq F^{*} \\
	&\bm{\hat x} \in \mathbb{R}^I_{+}  \\
	&	\bm{u},\bm{v},\bm{w}  \in \{0,1\} \\
	&\eqref{eq:primm}-\eqref{eq:primm_l} \quad \text{primal}\\
	&\eqref{eq:stattf}-\eqref{eq:stattf3} \quad \text{dual}\\
	&\eqref{eq:M1}-\eqref{eq:M2} \quad \text{slacks}\\
	&\eqref{eq:M4} \quad \text{frontier}.
\end{align*}
\endgroup

So far, we have not been very specific about the objective function $C(\bm{x}^0,\hat {\bm{x}})$.  Different functional forms can be introduced, and this may require the introduction of further variables to implement these. 

In Section \ref{sec:setting}, we approximated the firm's cost-of-change using combinations of the $\ell_0$ norm, the $\ell_1$ norm, and the squared $\ell_2$ norm. They are widely used in machine learning when close counterfactuals are sought in attempt to understand how to getter a more attractive outcome \citep{CARRIZOSAGroup2024}. The $\ell_0$ ``norm'', which strictly speaking is not a norm in the mathematical sense, counts the number of dimensions that has to be changed. The $\ell_1$ norm is the absolute value of the deviations. Lastly, $\ell_2^2$ is the Euclidean norm, that squares the deviations.

As a starting point, we therefore propose the following objective function:
\begin{equation}
	\label{eq:costf}
	C(\bm{x}^0,\hat {\bm{x}})=\nu_0 \|\bm{x}^0-\hat {\bm{x}}\|_0 +\nu_1 \|\bm{x}^0-\hat{ \bm{x}}\|_1 +\nu_2 \|\bm{x}^0-\hat{ \bm{x}}\|_2^2,
\end{equation}
where $\nu_0, \nu_1, \nu_2 \ge 0$.  Taking into account that there may be specific product input prices and output prices or that inputs may have varying degrees of difficulty to be changed, one can consider giving different weights to the deviations in each of the inputs.\\

In order to have a smooth expression of objective function \eqref{eq:costf}, additional decision variables and constraints have to be added to the counterfactual explanation problem in DEA. To linearize the $\ell_0$ norm, binary decision variables $\xi_i$ are introduced. For input $i$, $\xi_i=1$ models $x_i^0 \neq \hat x_i$, $i=1,\dots,I$. Using the big-M method the following constraints are added to our formulation:
\begin{align}
	-M_{\text{zero}}\xi_i \leq x_i^0&-\hat{x}_i \leq M_{\text{zero}} \xi_i, \quad   i=1,\dots,I \label{eq:l0}\\
	\xi_i &\in \{0,1\}, \quad  i=1,\dots,I, \label{eq:l0v}
\end{align}
where $M_{\text{zero}}$ is a sufficiently large constant.

For the $\ell_1$ norm we introduce continuous decision variables $\eta_i \ge 0$, $i=1,\dots,I$, to measure the absolute values of the deviations, $\eta_i=\left| x^0_i-\hat {x}_i\right|$, which is naturally implemented by the following constraints:
\begingroup
\allowdisplaybreaks
\begin{align}
	\eta_i \geq x_i^0-\hat{x}_i, \quad   i=1,\dots,I \label{eq:l1}\\
	-\eta_i \leq x_i^0-\hat{x}_i, \quad   i=1,\dots,I  \label{eq:l1v}\\
	\eta_i \geq 0, \quad   i=1,\dots,I.  \label{eq:l1nonneg}
\end{align}
\endgroup

Thus, the counterfactual explanation problem in DEA with cost function $C$ in \eqref{eq:costf}, hereafter \eqref{eq:cedeafinal}, reads as follows:
\begin{align}
	\label{eq:cedeafinal} \tag{CEDEA}
	\min_{\bm{\hat x},F, \bm{\beta},\bm{\gamma}_I,\bm{\gamma}_O, \bm{u},\bm{v},\bm{w},\bm{\eta},\bm{\xi}} \quad & \nu_0 \sum_{i=1}^I \xi_i +\nu_1 \sum_{i=1}^I \eta_i + \nu_2 \sum_{i=1}^I \eta_i^2\\
	\text{s.t.} \quad &F\leq F^{*} \nonumber \\
	&\bm{\hat x} \in \mathbb{R}^I_{+} \nonumber  \\
	&	\bm{u},\bm{v},\bm{w}  \in \{0,1\} \nonumber \\
	&\eqref{eq:primm}-\eqref{eq:primm_l}  ,\eqref{eq:M1}-\eqref{eq:M4},\eqref{eq:stattf}-\eqref{eq:stattf3},\nonumber\\
	& \eqref{eq:l0}-\eqref{eq:l0v}, \eqref{eq:l1}-\eqref{eq:l1nonneg}.\nonumber
\end{align}
Notice that in Problem \eqref{eq:cedeafinal} we assumed $\mathcal{X}(\bm{x}^0)=\mathbb{R}^I_{+}$ as the feasible space for $\bm{\hat x}$. Other relevant constraints for the counterfactual inputs could easily be added, e.g., bounds or relative bounds on the inputs, or inputs that cannot be changed in the short run, say capital expenses, or that represent environmental conditions beyond the control of the DMU.

In the case where only the $\ell_0$ and $\ell_1$ norms are considered, i.e., $\nu_2=0$, the objective function as well as the constraints are linear, while we have both binary and continuous decision variables. Therefore, Problem \eqref{eq:cedeafinal} can be solved using an Mixed Integer Linear Programming (MILP) solver. Otherwise, when $\nu_2\neq 0$,  Problem \eqref{eq:cedeafinal} is a Mixed Integer Convex Quadratic model with linear constraints, which can be solved with standard optimization packages. When all three norms are used, Problem \eqref{eq:cedeafinal} has $3I+K+2+O$ continuous variables and $2I+O+K+1$ binary decision variables. It has $7I+3O+3K+5$ constraints, plus the non-negativity and binary nature of the variables. The computational experiments show that this problem can be solved efficiently for our real-world dataset.

We can think of the objective function $C$ in different ways.

One possibility is to see it as an instrument to explore the production possibilities. The use of a combinations of the $\ell_0$, $\ell_1$  and $\ell_2$ norms seems natural here. Possible extensions could involve other $\ell_p$ norms,  $ \| \bm{x}^0-\hat {\bm{x}}\|_{p}:=\left(\sum _{i=1}^{I} \left| x^0_i-\hat{x}_i\right|^{p}\right)^{1/p}$. For all $p \in [1,\infty)$, $\ell_p$ is convex. This makes the use of $\ell_p$ norms convenient in generalizations of Problem \eqref{eq:cedeafinal}. Of course, arbitrary $\ell_p$ norms may lead to more complicated implementations in existing softwares since the objective function may no longer be quadratic.

Closely related to the instrumental view of the objective function is the idea of approximations. At least as a reasonable initial approximation of more complicated functions, many objective functions $C$ can be approximated by the form in \eqref{eq:costf}.

To end, one can link the form of $C$ closer to economic theory. In the economic literature there have been many studies on \emph{factor adjustments costs}. It is commonly believed that firms change their demand for inputs only gradually and with some delay, cf.\ 
e.g.~\cite{HamermestPfani1995}. For labor inputs, the factor adjustment costs include disruptions to production occurring when changing employment causes workers' assignments to be rearranged. Laying off or hiring new workers is also costly. There are search costs (advertising, screening, and processing new employees); the cost of training (including disruptions to production as previously trained workers' time is devoted to on-the-job instruction of new workers); severance pay (mandated and otherwise); and the overhead cost of maintaining that part of the personnel function dealing with recruitment and worker outflows. Following again \cite{HamermestPfani1995}, the literature on both labor and capital goods adjustments has overwhelmingly relied on one form of $C$, namely that of symmetric convex adjustment costs much like we use in \eqref{eq:costf}. Indeed, in the case of only one production factor, the most widely used function form is simply the quadratic one. Hall \cite{Hall2004} and several others have  tried to estimate the costs of adjusting labor and capital inputs.  Using a Cobb-Douglas production function, and absent adjustment costs and absent changes in the ratios of factor prices, an increase in demand or in another determinant of industry equilibrium would cause factor inputs to change in the same proportion as outputs. Adjustment costs are introduced as reductions in outputs and are assumed to depend on the squared growth rates in labor and capital inputs - the larger the percentage change, the larger the adjustment costs. Another economic approach to the cost-of-change modelling is to think of \emph{habits}. In firms - as in the private life - habits are useful. In the performance of many tasks, including complicated ones, it is easiest to go into automatic mode and let a behavior unfold. When an efficiency requirement is introduced, habits may need to change and this is costly. The relevance and strength of habit formation has also been studied empirically using panel data, cf.\ e.g.~\cite{Dynan2000} and the references herein. Habit formation should ideally be considered in a dynamic framework. To keep it simple, we might consider two periods - the past, where $\bm{x}^0$ was used and the present, where $\hat{\bm x}$ is consumed.  The utility in period two will then typically depend on the difference or ratio of present to past consumption, $\hat{\bm x}-\bm{x}^0$ or, in the unidimensional case, $\hat{\bm{x}} / \bm{x}^0$. Examples of functional forms one can use are provided in for example \cite{Fuhrer2000}.

\section{A banking application}
\label{sec:num}

In this section, we illustrate our methodology using real-world data on bank branches, \cite{SchaffnitEA1997}, by constructing a collection of counterfactual explanations for each of the inefficient firms that can help them learn about the DEA benchmarking model and how they can improve their efficiency.

The data is described in more detail in Section \ref{sec:data}, where we consider a model of bank branch production with $I=5$ inputs and $O=3$ outputs, and thus a production possibility set in $\mathbb{R}^8_{+}$, spanned by $K+1=267$ firms. In Section \ref{sec:illustrations}, we will focus on changing the inputs, and therefore the counterfactual explanations will be obtained with Problem \eqref{eq:cedeafinal}. We will discuss the results obtained with different cost of change functions $C$, reflecting the effort an inefficient firm will need to spend to change to its counterfactual instance, and different desired levels of efficiency $E^*$. The Farrell projection discussed in Section \ref{sec:setting} is added for reference. The counterfactual analysis sheds light on the nature of the DEA benchmarking model, which is otherwise hard to comprehend because of the many firms and inputs and outputs involved in the construction of the technology.

All optimization models have been implemented using Python 3.8 and as solver Gurobi 9.0 \citep{gurobi}. We have solved Problem \eqref{eq:cedeafinal} with $M_I= M_0=M_f=1000$ and $M_{\text{zero}}=1$. The validity of $M_{\text{zero}}=1$ will be shown below. Our numerical experiments have been conducted on a PC, with an Intel R CoreTM i7-1065G7 CPU @ 1.30GHz 1.50 GHz processor and 16 gigabytes RAM. The operating system is 64 bits.

\subsection{The data}
\label{sec:data}

The data consist of five staff categories and three different types of outputs in the Ontario branches of a large Canadian bank. The inputs are measured as full-time equivalents (FTEs), and the outputs are the average monthly counts of the different transactions and maintenance activities. Observations with input values equal to 0 are removed, leaving us with an actual dataset with 267 branches. Summary statistics are provided in Table~\ref{Tab:SummaryStats}.

\begin{table}[h]
	\footnotesize
	\begin{center}	
		\begin{tabular}{lrrrrr}
			\hline
			&	Mean	&Min	&Max	&Std. dev.  \\ \hline	
			INPUTS & & &   \\ \hline
			Teller &	5.83&	0.49&	39.74&	3.80 \\
			Typing&	1.05	&0.03&	22.92&	1.84\\
			Accounting \& ledgers	&4.69&	0.80&	65.93&	5.13\\
			Supervision&	2.05&	0.43&	38.29	&2.66	\\
			Credit&	4.40&	0.35&	55.73&	6.19\\ \hline
			OUTPUTS& & & & \\ \hline
			Term accounts&	2788&	336&	22910	&2222\\
			Personal loan accounts&	117&	0&	1192&	251\\
			Commercial loan accounts&	858&	104&	8689&	784 \\
			\hline
		\end{tabular}
		\caption{Descriptive statistics of the Canadian bank branches dataset in \cite{SchaffnitEA1997}}\label{Tab:SummaryStats}
	\end{center}	
\end{table}

After calculating all the efficiencies through Problem \eqref{eq:primal}, one has that 236 firms of the 267 ones are inefficient. Out of those, 219 firms have an efficiency below 90\%, 186 below 80\%, 144 below 70\%, 89 below 60\% and 49 below 50\%.

\subsection{Counterfactual analysis of bank branches}
\label{sec:illustrations}

To examine the inefficient firms, we will determine counterfactual explanations for these. Prior to that, we have divided each input by its maximum value across all firms. We notice that this has no impact on the solution since DEA models are invariant to linear transformations of inputs and outputs. Also, this makes valid choosing $M_{\text{zero}}=1$, since the values of all inputs are upper bounded by 1.

\begin{figure}
	\centering
	\includegraphics[width=0.575\linewidth]{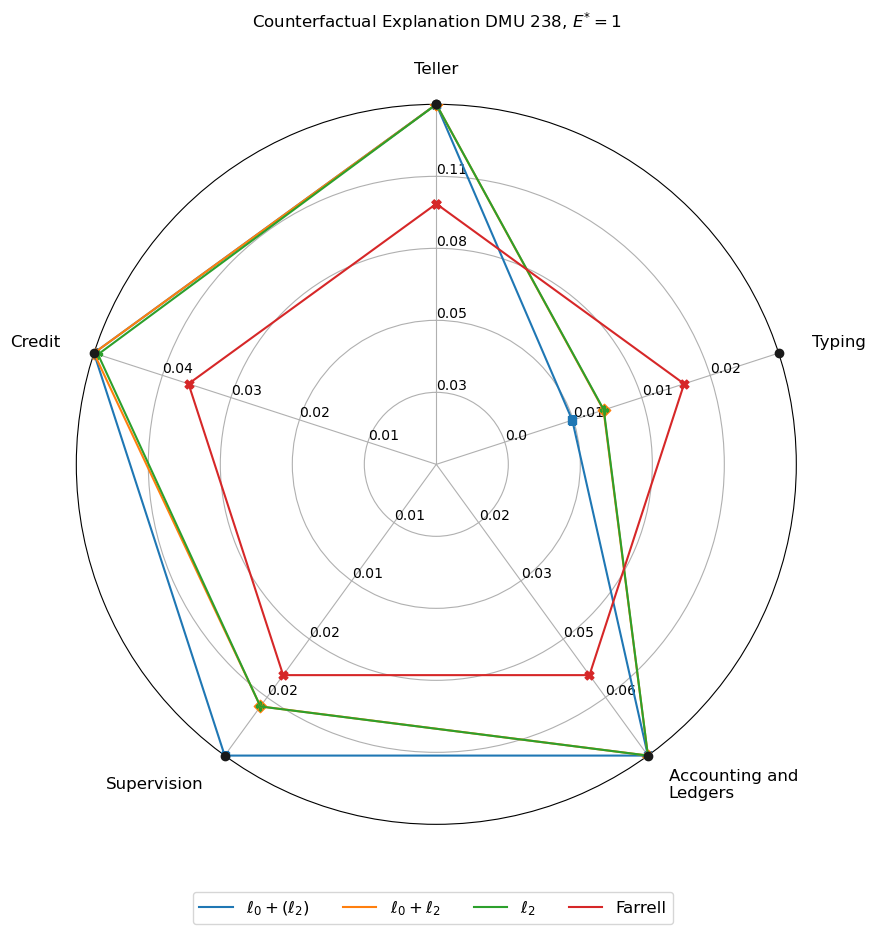}
	\caption{Counterfactual Explanations for firm 238 with Problem \eqref{eq:cedeafinal} and desired efficiency $E^{*}=1$.}
	\label{fig:spider1}
\end{figure}

\begin{figure}
	\subfloat[$C=\ell_0+(\ell_2)$\label{fig:heatmapl0}]
	{\includegraphics[width=.20\linewidth]{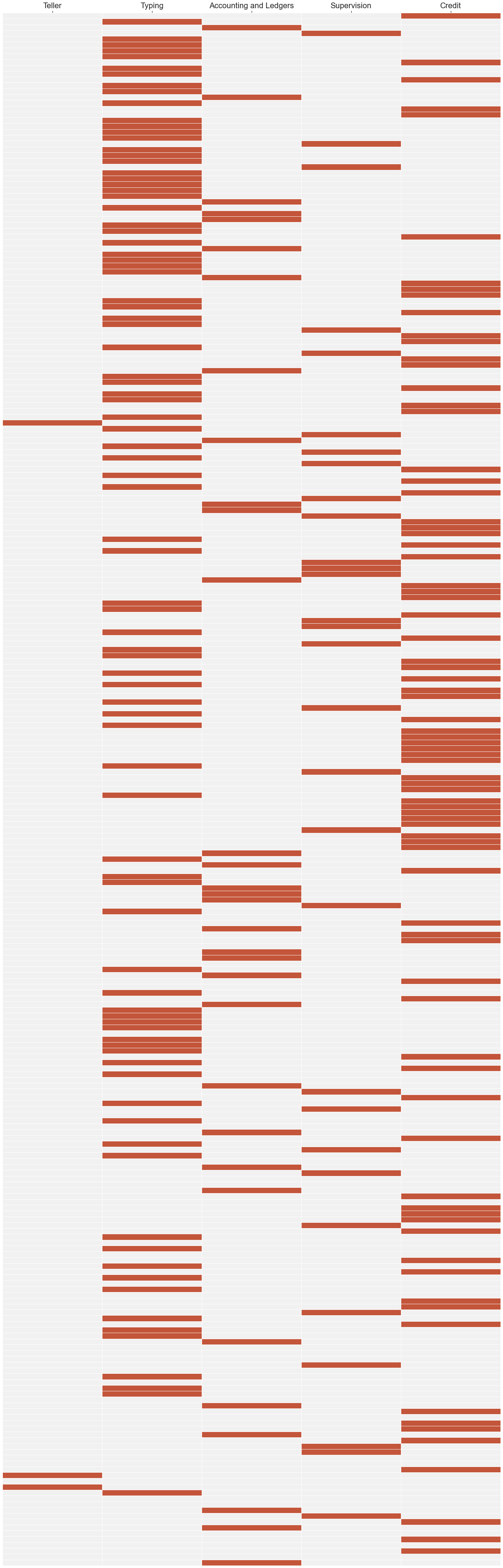}}\hfill
	\subfloat[$C=\ell_0+\ell_2$\label{fig:heatmapl0l2}]
	{\includegraphics[width=.20\linewidth]{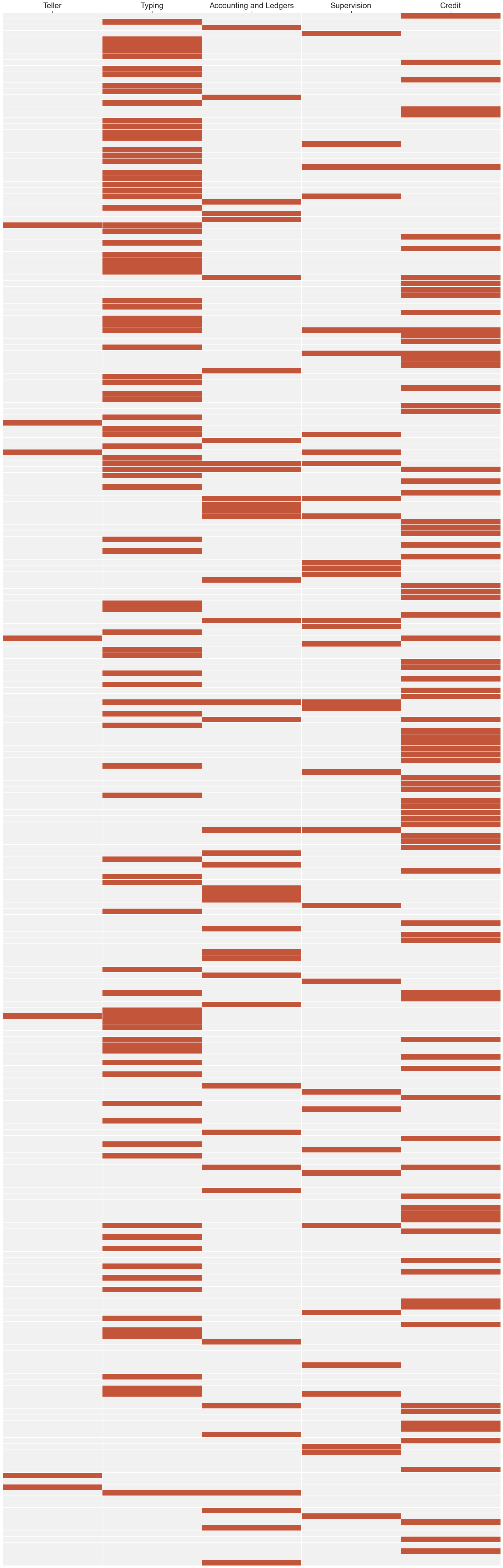}}\hfill
	\subfloat[$C=\ell_2$\label{fig:heatmapl2}]
	{\includegraphics[width=.20\linewidth]{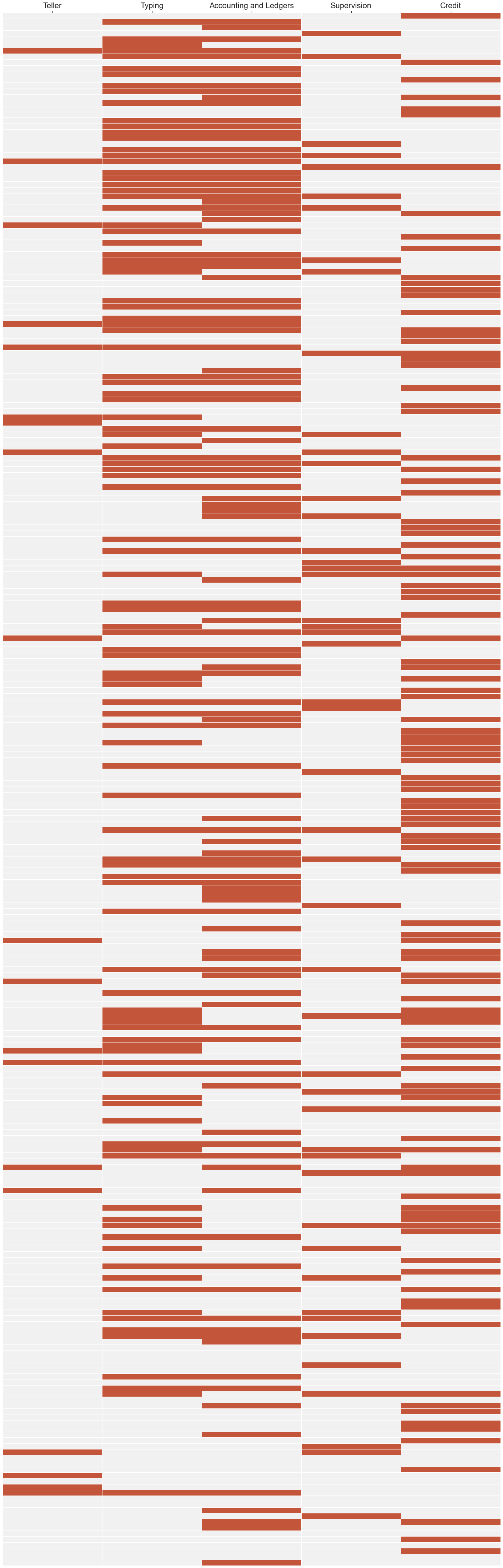}}
	\caption{The inputs that change when we impose a desired efficiency of $E^{*}=1$}
	\label{fig:heatmap}
\end{figure}

\begin{figure}
	\centering
	\includegraphics[width=0.575\linewidth]{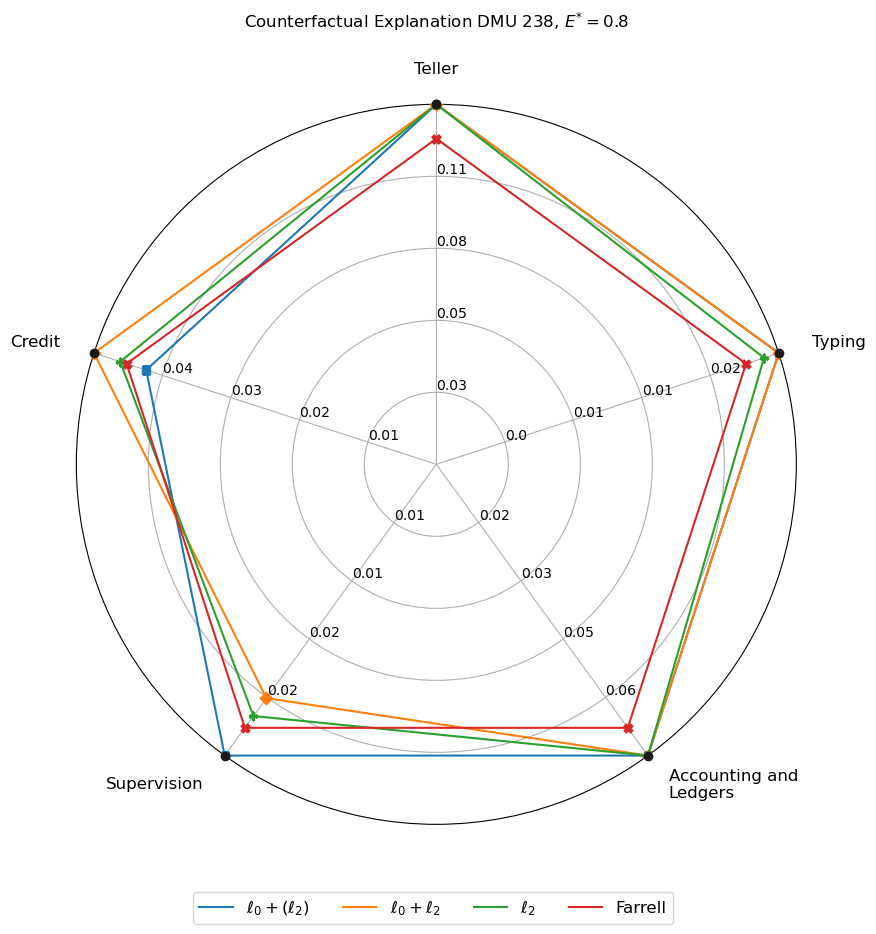}
	\caption{Counterfactual Explanations for DMU 238 with Problem \eqref{eq:cedeafinal} and desired efficiency $E^{*}=0.8$}
	\label{fig:spider08}
\end{figure}

\begin{figure}
	\subfloat[$C=\ell_0+(\ell_2)$\label{fig:heatmapl0_2}]
	{\includegraphics[width=.20\linewidth]{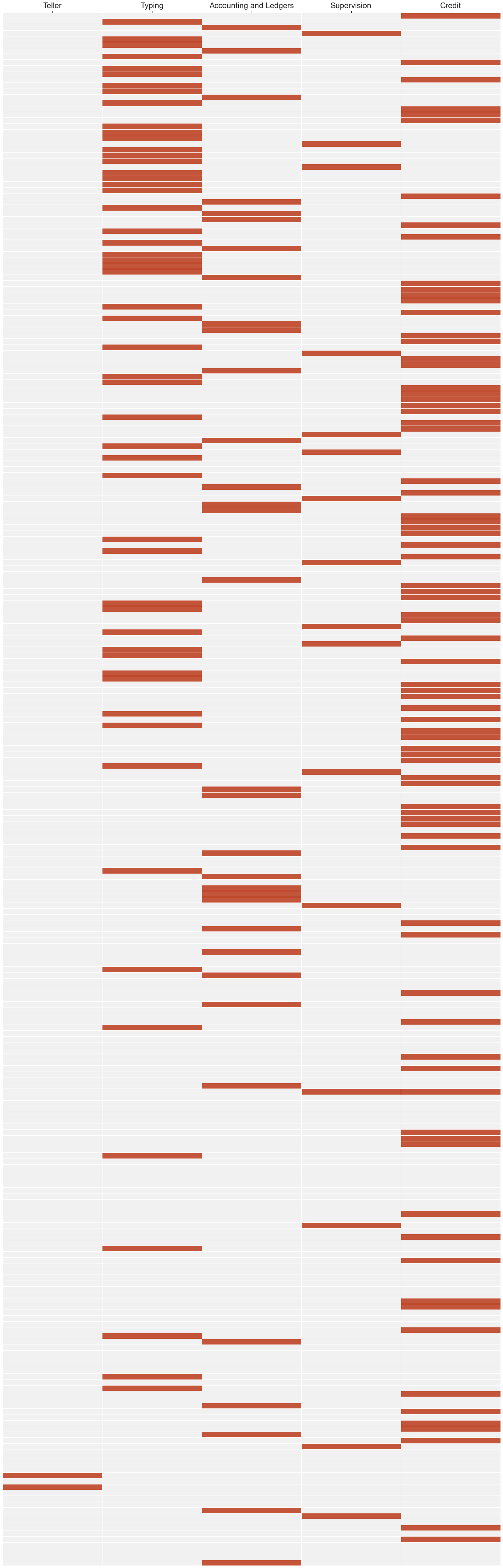}}\hfill
	\subfloat[$C=\ell_0+\ell_2$\label{fig:heatmapl0l2_2}]
	{\includegraphics[width=.20\linewidth]{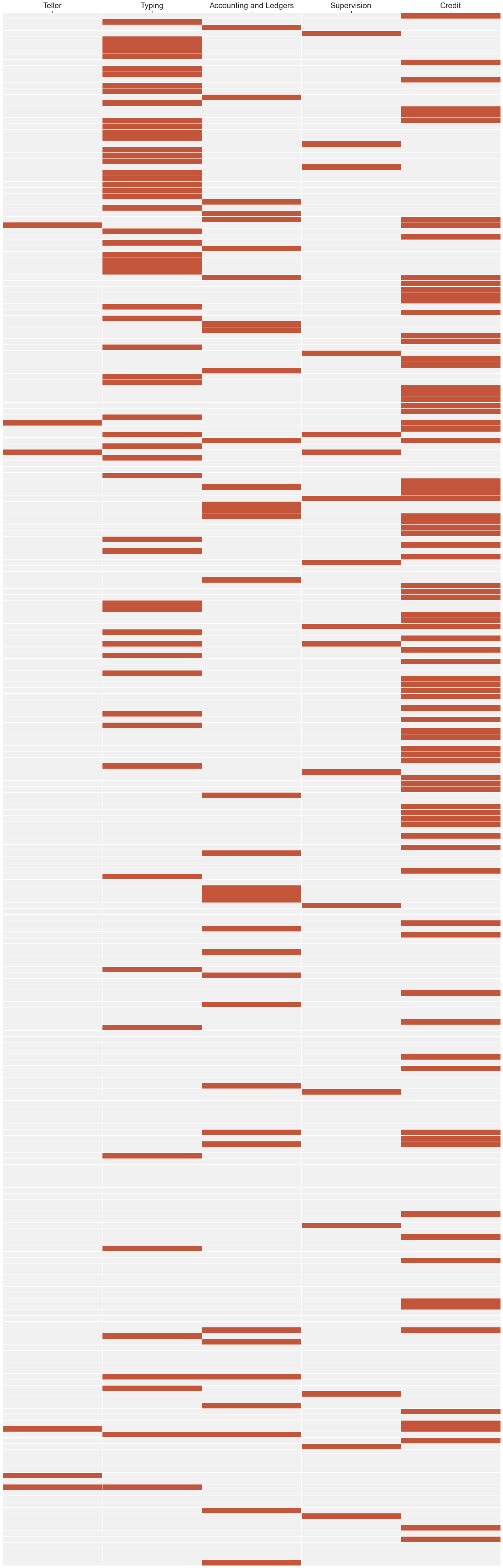}}\hfill
	\subfloat[$C=\ell_2$\label{fig:heatmapl2_2}]
	{\includegraphics[width=.20\linewidth]{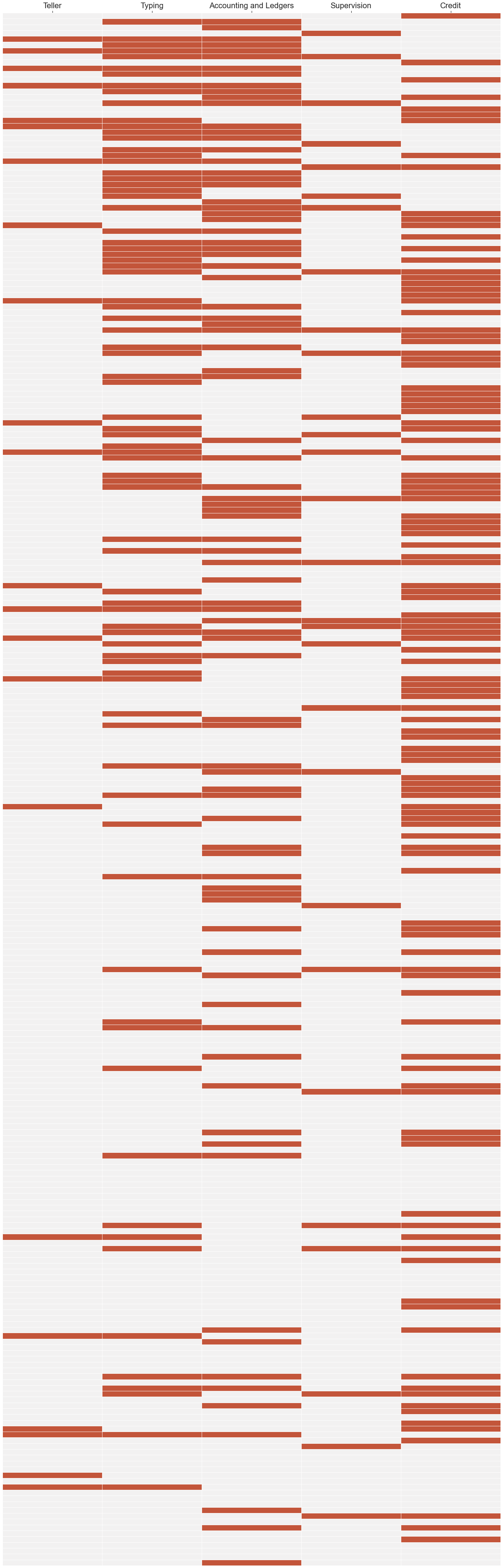}}
	\caption{The inputs that change when we impose a desired efficiency of $E^{*}=0.8$}
	\label{fig:heatmap2}
\end{figure}

We will use three different cost functions, by changing the values of the parameters $\nu_0,\nu_1,\nu_2$ in \eqref{eq:costf}, as well as two different values of the desired efficiency of the counterfactual instance, namely $E^{*}=1$ and $0.8$. In the first implementation of the cost function, which we denote $\ell_0 +(\ell_2)$, we use $\nu_0=1$, $\nu_2=10^{-3}$ and $\nu_1=0$, i.e., we will seek to minimize the $\ell_0$ norm and only introduce a little bit of the squared Euclidean norm to ensure a unique solution of Problem \eqref{eq:cedeafinal}. In the second implementation, which we call $\ell_0+\ell_2$, we take $\nu_0=1$, $\nu_1=0$ and $\nu_2=10^5$, such that the squared Euclidean norm has a higher weight than in cost function $\ell_0 +(\ell_2)$. Finally, we denote by $\ell_2$ the cost function that focuses on the minimization of the squared Euclidean norm only, i.e., $\nu_0=\nu_1=0$ and $\nu_2=1$. The summary of all the cost functions used can be seen in Table \ref{tab:costf}. Calculations were also done for the $\ell_1$ norm, i.e., $\nu_0=\nu_2=0$ and $\nu_1=1$, but as the solutions found were similar to those for cost function $\ell_0+\ell_2$, for the sake of clarity of presentation, they are omitted. We start the discussion of the counterfactual explanations obtained with $E^{*}=1$, as summarized in Figures \ref{fig:spider1}-\ref{fig:heatmap} and Tables \ref{tab:changes1}-\ref{tab:changes1_2}. We then move on to a less demanding desired efficiency, namely, $E^{*}=0.8$. These results are summarized in Figures \ref{fig:spider08}-\ref{fig:heatmap2} and Tables \ref{tab:changes2}-\ref{tab:changes2_2}.

\begin{table}[h]
	\footnotesize
	\begin{center}	
		\begin{tabular}{l|ccc}
			Cost function & $\nu_0$ & $\nu_1$ &$\nu_2$ \\ \hline
			\rule{0pt}{3ex}$\ell_0+(\ell_2)$& 1&0 &$ 10^{-3} $\\
			$\ell_0+\ell_2$ &1&0 &$10^5$\\
			$\ell_2$  &0&0&1
		\end{tabular}
		\caption{Value of the parameters $\nu_0,\nu_1$ and $\nu_2$ in \eqref{eq:costf} for the different cost functions used}
		\label{tab:costf}
	\end{center}	
\end{table}
Let us first visualize the counterfactual explanations for a specific firm. Consider, for instance, firm 238, which has an original efficiency of $E^{0}=0.72$. We can visualize the different counterfactual explanations generated by the different cost functions using  a spider chart, see Figure \ref{fig:spider1}. In addition to the counterfactual explanations obtained with Problem \eqref{eq:cedeafinal}, we also illustrate the so-called Farrell projection. In the spider chart, each axis represents an input and the original values of the firm corresponds to the outer circle. Figure \ref{fig:spider1} shows the different changes needed depending on the cost function used. With the $\ell_0+(\ell_2)$, where the focus is to mainly penalize the number of inputs changed, we see that only the typing personnel has to be changed, leaving the rest of the inputs unchanged. Nevertheless, because only one input is changed, it has to be decreased by 60\% from the original value. The Farrell solution decreases the typing personnel by 28\% of its value, but to compensate, it changes the remaining four inputs proportionally. When the $\ell_0+\ell_2$ cost function is used, the typing personnel keeps on needing to be changed, but the change is smaller, this time by 51\% of its value. The supervision personnel needs also to be decreased by 16\% of its value, while the rest of the inputs remain untouched. Increasing the weight on the Euclidean norm in the cost function gives us the combination of the two inputs that are crucial to change in order to gain efficiency, as well as the exact amount that they need to be reduced.  Finally, using only the Euclidean norm, the typing, supervision and credit personnel are the inputs to be changed, the typing input is reduced slightly less than with the $\ell_0+\ell_2$ in exchange of reducing just by 1\% the credit input. Notice that the teller and accounting and ledgers personnel are never changed in the counterfactual explanations generated by our methodology, which leads us to think that these inputs are not the ones leading firm 238 to have its original low efficiency.

The analysis above is for a single firm. We now present some statistics about the counterfactual explanations obtained for all the inefficient firms. Recall that these are 236 firms, and that Problem \eqref{eq:cedeafinal} has be solved for each of them. In Table \ref{tab:changes1} we show for each cost function, how often an input has to be changed. For instance, the value 0.09 in the last row of the Teller column shows that in 9\% of all firms, we have to change the number of tellers when the aim is to find a counterfactual instance using the Euclidean norm. When more weight is given to the $\ell_0$ norm, few inputs are changed. Indeed, for the Teller column, with the $\ell_0+\ell_2$, 3\% of all firms change it, instead of 9\%, and this number decreases to 1\% when the $\ell_0+(\ell_2)$ is used.  The same pattern can be observed in all inputs, particularly notable in the Acc.\ and Ledgers personnel, that goes from changing in more than half of the banks with the Euclidean norm, to changing in only 14\% of the firms. The last column of Table \ref{tab:changes1}, Mean $\ell_0(\bf{x}-\bf{\hat{x}})$, shows how many inputs are changed on average when we use the different cost functions. With the $\ell_0+(\ell_2)$ only one input has to be decreased, thus with this cost function one detects the crucial input to be modified to be fully efficient, leaving the rest fixed. In general, the results show that, for the inefficient firms, the most common changes leading to full efficiency is to reduce the number of typists and the number of credit officers. The excess of typists is likely related to the institutional setting. Bank branches need the so-called typists for judicial regulations, but they only need the services to a limited degree, see also Table \ref{Tab:SummaryStats}. In such cases, it may be difficult to match the full time equivalents employed precisely to the need. The excess of Credit officers is more surprising since, in particular, they are one of the best paid personnel groups.

In Table \ref{tab:changes1_2}, we look at the size of the changes and not just if a change has to take place or not. The interpretation of the value 0.43 under the first row and the Teller column suggests that when the teller numbers have to be changed, they are reduced by 43\% from the initial value, on average. Since several inputs may have to change simultaneously, defining the vector of the relative changes $\bm{r}=((x_i-\hat{x}_i)/x_i)_{i=1}^I$, the last column shows the mean value of the Euclidean norm of this vector. We see, for example, that in the relatively few cases the teller personnel has to change under $\ell_0+(\ell_2)$, the changes are relatively large. We see again the difficulties the bank branches apparently have hiring the right amount of typists. We saw in Table \ref{tab:changes1} that they often have to change and we see now that the changes are non-trivial with about a half excess full time equivalents.
\begin{table}[h]
	\footnotesize
	\begin{center}	
		\begin{tabular}{l|cccccc}
			Cost function & Teller & Typing &Acc.\ and Ledgers&Supervision&Credit
			& Mean $\ell_0(\bf{x}-\bf{\hat{x}})$ \\ \hline
			$\ell_0+(\ell_2)$& 0.01&0.38 &0.14 &0.13  & 0.34& 1.00\\
			$\ell_0+\ell_2$ &0.03&0.40 &0.17& 0.14&0.38&1.13\\
			$\ell_2$  & 0.09&0.45&0.51&0.21&0.47&1.72
		\end{tabular}
		\caption{Average results on how often the inputs (personnels) change when desired efficiency is $E^{*}=1$.}
		\label{tab:changes1}
	\end{center}	
\end{table}

\begin{table}[h]
	\footnotesize
	\begin{center}	
		\begin{tabular}{l|cccccc}
			Cost function & Teller & Typing &Acc.\ and Ledgers&Supervision&Credit
			& Mean $\ell_2(\bm{r})$ \\ \hline
			$\ell_0+(\ell_2)$& 0.43&0.61 &0.41 &0.43  & 0.37& 0.4743\\
			$\ell_0+\ell_2$ &0.21&0.58 &0.33& 0.38&0.35&0.4742\\
			$\ell_2$  & 0.11&0.53&0.14&0.27&0.29&0.4701
		\end{tabular}
		\caption{Average results on how much the inputs (personnels) change when desired efficiency is $E^{*}=1$.}
		\label{tab:changes1_2}
	\end{center}
\end{table}

In Figure \ref{fig:heatmap}, we use a heatmap to illustrate which input factors have to change for the individual firms using the three different cost functions in Table \ref{tab:changes1}. Rows with no markings represent firms that were fully efficient to begin with. We see as we would expect that the more weight we put on the Euclidean norm, the more densely populated the illustration becomes, i.e., the more inputs have to change simultaneously.

So far we have asked for counterfactual instances that are fully efficient. If we instead only ask for the counterfactual instances to be at least 80\% efficient, only 186 firms need to be studied.
As before, let us first visualize the counterfactual explanations for firm 238, which had an original efficiency of $E^{0}=0.72$. In Figure \ref{fig:spider08}, we can see the different changes when imposing $E^{*}=0.8$. We again see that the Farrell approach reduces all inputs proportionally, specifically by 9.5\% of their values. We see also that under the $\ell_0+(\ell_2)$ norm, only Credit personnel has to be reduced, by 15\%. Under the Euclidean norm, Teller and Acc.\ and Ledgers personnel are not affected while Typing, Supervision and Credit officers have to be saved, by 4\%, 13\% and 7\%, respectively. Notice that only the change in Supervision is higher in this case than in the Farrell solution, while the decrease in the remain four inputs is significantly smaller for the Euclidean norm. Recall that in Figure \ref{fig:spider1} the counterfactual explanations for the same firm 238 have been calculated imposing $E^{*}=1$. Altering the desired efficiency level from $E^*=0.8$ to $E^*=1$ leads to rather dramatic changes in the counterfactual explanations. For the $\ell_0+(\ell_2)$ cost function, for a desired efficiency of $E^*=0.8$, we needed to decrease the Credit personnel dramatically whereas for a desired efficiency of $E^*=1$, it is suggested to leave unchanged the Credit personnel and to change the Typing personnel instead. On the other hand, what remains the same is the fact that Teller and Acc.\ and Ledgers officers are never affected in the counterfactual explanations with the three cost functions.

After the analysis for a single firm, now we present statistics about the counterfactual explanations obtained for all 168 firms that had an original efficiency below 80\%. The frequency of changes and the relative sizes of the changes are shown in Tables \ref{tab:changes2} and \ref{tab:changes2_2}. We see, as we would expect, that the amount of changes necessary is reduced. On the other hand, the inputs to be changed are not vastly different. The tendency to change in particular credit officers is slightly larger now.

\begin{table}[h]
	\footnotesize
	\begin{center}	
		\begin{tabular}{l|cccccc}
			Cost function & Teller & Typing &Acc.\ and Ledgers &Supervision &Credit & Mean $\ell_0(\bf{x}-\bf{\hat{x}})$ \\ \hline
			$\ell_0+(\ell_2)$& 0.01	&0.30&	0.18&	0.08&	0.44&	1.00
			\\
			$\ell_0+\ell_2$ &0.03&	0.32&	0.19&	0.09&	0.47&	1.11
			\\
			$\ell_2$  & 0.13&	0.43&	0.51&	0.18&	0.63&	1.88
			
		\end{tabular}
		\caption{Average results on how often the inputs (personnels) change when desired efficiency is $E^{*}=0.8$}
		\label{tab:changes2}
	\end{center}
\end{table}

\begin{table}[h]
	\footnotesize
	\begin{center}	
		\begin{tabular}{l|cccccc}
			Cost function & Teller & Typing &Acc.\ and Ledgers&Supervision&Credit
			& Mean $\ell_2(\bm{r})$ \\ \hline
			$\ell_0+(\ell_2)$& 0.28&	0.56&	0.28&	0.41&	0.27&	0,3708
			\\
			$\ell_0+\ell_2$ &0.12&	0.52&	0.25&	0.39&	0.26&	0.3707
			\\
			$\ell_2$  & 0.05&	0.41&	0.11&	0.25&	0.21&	0.3702
			
		\end{tabular}
		\caption{Average results on how much the inputs (personnels) change when desired efficiency is $E^{*}=0.8$}
		\label{tab:changes2_2}
	\end{center}
\end{table}

In Figure \ref{fig:heatmap2}, we show the input factors that need to change for the individual firms using the three different cost functions in Table \ref{tab:changes2} for the case now with $E^{*}=0.8$. As expected, we can see now an increasing number of rows with no markings compared to Figure \ref{fig:heatmap}, belonging to the firms that had already an efficiency of 0.8.

\section{Conclusions}
\label{sec:conclusion}

In this paper, we have proposed a collection of optimization models to setting targets and finding counterfactual explanations in DEA models, i.e., the least costly changes in the inputs or outputs of a firm that leads to a pre-specified (higher) efficiency level. With our methodology, we are able to include different ways to measure the proximity between a firm and its counterfactual, namely, using the  $\ell_0$,  $\ell_1$, and $\ell_2$ norms or a combination of them. Calculating counterfactual explanations involves finding ``close" alternatives in the complement of a convex set. We have reformulated this bilevel optimization problem as either an MILP or a Mixed Integer Convex Quadratic Problem with linear constraints. In our numerical section, we can see that for our banking application, we are able to solve this model to optimality.

DEA models can capture very complex relationships between multiple inputs and outputs. This allows more substantial evaluations and also offers a framework that can support many operational, tactical and strategic planning efforts. However, there is also a risk that such a model is seen as a pure black box which in turn can lead to mistrust and some degree of model or algorithm aversion. By looking at counterfactuals, a firm can get a better understanding of the production space and is more likely to trust in the modelling.

Counterfactuals in DEA can also help a firm choose which changes to implement. It is not always enough to simply think of a strategy and which factors can easily be changed, say a direction in input space. It is also important how the technology looks like and therefore how large such changes need to be to get a desired improvement in efficiency. In this way, the analysis of close counterfactuals can help endogenize the choice of both desirable and effective directions to move in. By varying the parameters of the cost function, the firm can even get a menu of counterfactuals, from which it can choose, having thus more flexibility and leading the evaluated firm to gain more trust in the underlying model.

Note also that by calculating the counterfactual explanations for all firms involved, as we did in our banking application, one can determine which combinations of inputs and outputs that most commonly shall be changed to improve efficiency. This is interesting from an overall system point of view. Society at large - or for example a regulator tasked primarily with incentivizing natural monopolies to improve efficiency - may not solely be interested that everyone becomes efficient. It may as well be important how the efficiency is improved, e.g.\ by reducing the use of imported or domestic resources or by laying off some particular types of labor and not other types.

There are several interesting extensions that can be explored in future research. Here we just mention two. One possibility is to use alternative efficiency measures to constrain the search for counterfactual instances. We have here used Farrell efficiency, which is by far the most common efficiency measure in DEA studies, but one might consider other alternative measures, e.g.\ additive ones like the excess measure. Another relevant extension could be to make the counterfactuals less individualized. One could for example look for the common features that counterfactual explanations should change across all individual firms and that lead to the minimum total cost.

\section*{Acknowledgements} This research has been financed in part by research projects EC H2020 MSCA RISE NeEDS (Grant 822214); COST Action CA19130 - FinAI; FQM-329, P18-FR-2369 and US-1381178 (Junta de Andaluc\'{\i}a), PID2019-110886RB-I00 and PID2022-137818OB-I00 (Ministerio de Ciencia, Innovaci\'on y Universidades, Spain), and Independent Research Fund Denmark (Grant 9038-00042A) - ``Benchmarking-based incentives and regulatory applications''. The support is gratefully acknowledged. 


\bibliographystyle{plainnat} 
\bibliography{bibliografia}	

\begin{thebibliography}{47}
\providecommand{\natexlab}[1]{#1}
\providecommand{\url}[1]{\texttt{#1}}
\expandafter\ifx\csname urlstyle\endcsname\relax
  \providecommand{\doi}[1]{doi: #1}\else
  \providecommand{\doi}{doi: \begingroup \urlstyle{rm}\Url}\fi

\bibitem[Agrell and Bogetoft(2017)]{AgrellBogetoftHANDBOOK}
P.J. Agrell and P.~Bogetoft.
\newblock Theory, techniques and applications of regulatory benchmarking and
  productivity analysis.
\newblock In \emph{Oxford Handbook of Productivity Analysis}, pages 523--555.
  Oxford University Press: Oxford, 2017.

\bibitem[Antle and Bogetoft(2019)]{AntleBogetoft2019}
R.~Antle and P.~Bogetoft.
\newblock Mix stickiness under asymmetric cost information.
\newblock \emph{Management Science}, 65\penalty0 (6):\penalty0 2787--2812,
  2019.

\bibitem[Aparicio and Pastor(2013)]{AparicioPastor2013}
J.~Aparicio and J.T. Pastor.
\newblock A well-defined efficiency measure for dealing with closest targets in
  {DEA}.
\newblock \emph{Applied Mathematics and Computation}, 219:\penalty0 9142--9154,
  2013.

\bibitem[Aparicio et~al.(2007)Aparicio, Ruiz, and
  Sirvent]{AparicioRuizSirvent2007}
J.~Aparicio, J.~L. Ruiz, and I.~Sirvent.
\newblock Closest targets and minimum distance to the pareto-efficient frontier
  in {DEA}.
\newblock \emph{Journal of Productivity Analysis}, 28:\penalty0 209--218, 2007.

\bibitem[Aparicio et~al.(2014)Aparicio, Mahlberg, Pastor, and
  Sahoo]{Aparicio_et_al.2014c}
J.~Aparicio, B.~Mahlberg, J.T. Pastor, and B.K. Sahoo.
\newblock Decomposing technical inefficiency using the principle of least
  action.
\newblock \emph{European Journal of Operational Research}, 239:\penalty0
  776--785, 2014.

\bibitem[Aparicio et~al.(2017)Aparicio, Cordero, and
  Pastor]{Aparicio_et_all_2017Omega}
J.~Aparicio, J.~M. Cordero, and J.~T. Pastor.
\newblock The determination of the least distance to the strongly efficient
  frontier in data envelopment analysis oriented models: modelling and
  computational aspects.
\newblock \emph{Omega}, 71:\penalty0 1--10, 2017.

\bibitem[Bogetoft(2012)]{Bogetoft2012}
P~Bogetoft.
\newblock \emph{Performance Benchmarking - Measuring and Managing Performance}.
\newblock Springer, New York, 2012.

\bibitem[Bogetoft and Hougaard(1999)]{Bogetoft1999}
P.~Bogetoft and J.L. Hougaard.
\newblock Efficiency evaluation based on potential (non-proportional)
  improvements.
\newblock \emph{Journal of Productivity Analysis}, 12:\penalty0 233--247, 1999.

\bibitem[Bogetoft and Otto(2011)]{BogetoftOtto2011}
P.~Bogetoft and L.~Otto.
\newblock \emph{{Benchmarking with DEA, SFA, and R}}.
\newblock Springer, New York, 2011.

\bibitem[Carrizosa et~al.(2024{\natexlab{a}})Carrizosa, {Ram\'{\i}rez Ayerbe},
  and {Romero Morales}]{carrizosaESWA24}
E.~Carrizosa, J.~{Ram\'{\i}rez Ayerbe}, and D.~{Romero Morales}.
\newblock Generating collective counterfactual explanations in score-based
  classification via mathematical optimization.
\newblock \emph{Expert Systems With Applications}, 238:\penalty0 121954,
  2024{\natexlab{a}}.

\bibitem[Carrizosa et~al.(2024{\natexlab{b}})Carrizosa, Ramírez-Ayerbe, and
  {Romero Morales}]{CARRIZOSAGroup2024}
E.~Carrizosa, J.~Ramírez-Ayerbe, and D.~{Romero Morales}.
\newblock Mathematical optimization modelling for group counterfactual
  explanations.
\newblock \emph{\emph{Forthcoming in} European Journal of Operational
  Research}, 2024{\natexlab{b}}.
\newblock URL \url{https://doi.org/10.1016/j.ejor.2024.01.002}.

\bibitem[Charnes et~al.(1978)Charnes, Cooper, and Rhodes]{Charnes1978}
A.~Charnes, William~W. Cooper, and E.~Rhodes.
\newblock Measuring the efficiency of decision making units.
\newblock \emph{European Journal of Operational Research}, 2:\penalty0
  429--444, 1978.

\bibitem[Charnes et~al.(1979)Charnes, Cooper, and Rhodes]{Charnes1979}
A.~Charnes, W.~W. Cooper, and E.~Rhodes.
\newblock {Short Communication: Measuring the Efficiency of Decision Making
  Units}.
\newblock \emph{European Journal of Operational Research}, 3:\penalty0 339,
  1979.

\bibitem[Charnes et~al.(1995)Charnes, Cooper, Lewin, and
  Seiford]{Coopereabook1995}
A.~Charnes, W.~W. Cooper, A.~Y. Lewin, and L.~M. Seiford.
\newblock \emph{Data Envelopment Analysis: Theory, Methodology and
  Applications}.
\newblock Kluwer Academic Publishers, Boston, USA, 1995.

\bibitem[Cherchye et~al.(2013)Cherchye, Rock, Dierynck, Roodhooft, and
  Sabbe]{Cherchye_ea2013}
L.~Cherchye, B.~De Rock, B.~Dierynck, F.~Roodhooft, and J.~Sabbe.
\newblock Opening the ``black box'' of efficiency measurement: Input allocation
  in multioutput settings.
\newblock \emph{Operations Research}, 61\penalty0 (5):\penalty0 1148--1165,
  2013.

\bibitem[Du et~al.(2019)Du, Liu, and Hu]{du2019techniques}
M.~Du, N.~Liu, and X.~Hu.
\newblock Techniques for interpretable machine learning.
\newblock \emph{Communications of the ACM}, 63\penalty0 (1):\penalty0 68--77,
  2019.

\bibitem[Dynan(2000)]{Dynan2000}
K.E. Dynan.
\newblock Habit formation in consumer preferences: Evidence from panel data.
\newblock \emph{American Economic Review}, 90\penalty0 (3):\penalty0 391--406,
  2000.

\bibitem[{European Commission}(2020)]{EUwhitepaperAI20}
{European Commission}.
\newblock \emph{{White Paper on Artificial Intelligence : a European approach
  to excellence and trust}}.
\newblock
  {\url{https://ec.europa.eu/info/publications/white-paper-artificial-intelligence-european-approach-excellence-and-trust_en}},
  2020.

\bibitem[F\"{a}re and Grosskopf(2000)]{FareGrosskopf2000}
R.~F\"{a}re and S.~Grosskopf.
\newblock Network {DEA}.
\newblock \emph{Socio-Economic Planning Sciences}, 34\penalty0 (1):\penalty0
  35--49, 2000.

\bibitem[F\"{a}re et~al.(2013)F\"{a}re, Grosskopf, and
  Whittaker]{FareGrosskopfWittaker2013}
R.~F\"{a}re, S.~Grosskopf, and G.~Whittaker.
\newblock Directional output distance functions: endogenous directions based on
  exogenous normalization constraints.
\newblock \emph{Journal of Productivity Analysis}, 40:\penalty0 267--269, 2013.

\bibitem[F\"{a}re et~al.(2017)F\"{a}re, Pasurkac, and M.Vardanyan]{Fare_ea2017}
R.~F\"{a}re, C.~Pasurkac, and M.Vardanyan.
\newblock On endogenizing direction vectors in parametric directional distance
  function-based models.
\newblock \emph{European Journal of Operational Research}, 262:\penalty0
  361--369, 2017.

\bibitem[Fischetti and Jo(2018)]{Fischetti2018}
M.~Fischetti and J.~Jo.
\newblock Deep neural networks and mixed integer linear optimization.
\newblock \emph{Constraints}, 23\penalty0 (3):\penalty0 296--309, 2018.

\bibitem[Fuhrer(2000)]{Fuhrer2000}
J.~C. Fuhrer.
\newblock Habit formation in consumption and its implications for
  monetary-policy models.
\newblock \emph{The American Economic Review}, 90\penalty0 (4):\penalty0
  367--390, 2000.

\bibitem[Goodman and Flaxman(2017)]{goodmanAIM17}
B.~Goodman and S.~Flaxman.
\newblock {European Union regulations on algorithmic decision-making and a
  ``right to explanation''}.
\newblock \emph{AI Magazine}, 38\penalty0 (3):\penalty0 50--57, 2017.

\bibitem[Guidotti(2022)]{guidotti2022counterfactual}
R.~Guidotti.
\newblock Counterfactual explanations and how to find them: literature review
  and benchmarking.
\newblock \emph{\emph{Forthcoming in} Data Mining and Knowledge Discovery},
  2022.

\bibitem[Gurobi~Optimization(2021)]{gurobi}
LLC Gurobi~Optimization.
\newblock Gurobi optimizer reference manual, 2021.
\newblock URL \url{http://www.gurobi.com}.

\bibitem[Hall(2004)]{Hall2004}
R.E. Hall.
\newblock Measuring factor adjustment costs.
\newblock \emph{The Quarterly Journal of Economics}, 119\penalty0 (3):\penalty0
  899--927, 2004.

\bibitem[Hamermesh and Pfann(1999)]{HamermestPfani1995}
D.~S. Hamermesh and G.~A. Pfann.
\newblock Adjustment costs in factor demand.
\newblock \emph{Journal of Economic Literature}, 34\penalty0 (3):\penalty0
  1264--1292, 1999.

\bibitem[Haney and Pollitt(2009)]{Haney2009}
A.B. Haney and M.G. Pollitt.
\newblock Efficiency analysis of energy networks: An international survey of
  regulators.
\newblock \emph{Energy Policy}, 37\penalty0 (12):\penalty0 5814--5830, 2009.

\bibitem[Kao(2009)]{Kao2009}
C.~Kao.
\newblock Efficiency decomposition in network data envelopment analysis: A
  relational model.
\newblock \emph{European Journal of Operational Research}, 192:\penalty0
  949--962, 2009.

\bibitem[Karimi et~al.(2022)Karimi, Barthe, Sch\"{o}lkopf, and
  Valera]{karimi2022survey}
A.-H. Karimi, G.~Barthe, B.~Sch\"{o}lkopf, and I.~Valera.
\newblock A survey of algorithmic recourse: contrastive explanations and
  consequential recommendations.
\newblock \emph{{ACM} Computing Surveys}, 55\penalty0 (5):\penalty0 1--29,
  2022.

\bibitem[Lundberg and Lee(2017)]{lundberg2017unified}
S.M. Lundberg and S.-I. Lee.
\newblock A unified approach to interpreting model predictions.
\newblock In \emph{Advances in Neural Information Processing Systems}, pages
  4765--4774, 2017.

\bibitem[Martens and Provost(2014)]{martensMISQ14}
D.~Martens and F.~Provost.
\newblock Explaining data-driven document classifications.
\newblock \emph{MIS Quarterly}, 38\penalty0 (1):\penalty0 73--99, 2014.

\bibitem[Molnar et~al.(2020)Molnar, Casalicchio, and
  Bischl]{molnar2020interpretable}
C.~Molnar, G.~Casalicchio, and B.~Bischl.
\newblock Interpretable machine learning--a brief history, state-of-the-art and
  challenges.
\newblock In \emph{Joint European Conference on Machine Learning and Knowledge
  Discovery in Databases}, pages 417--431. Springer, 2020.

\bibitem[Parmentier and Vidal(2021)]{parmentier2021optimal}
A.~Parmentier and T.~Vidal.
\newblock Optimal counterfactual explanations in tree ensembles.
\newblock In \emph{International Conference on Machine Learning}, pages
  8422--8431. PMLR, 2021.

\bibitem[Parmeter and Zelenyuk(2019)]{parmeter2019combining}
C.F. Parmeter and V.~Zelenyuk.
\newblock Combining the virtues of stochastic frontier and data envelopment
  analysis.
\newblock \emph{Operations Research}, 67\penalty0 (6):\penalty0 1628--1658,
  2019.

\bibitem[Petersen(2018)]{petersen2018directional}
N.C. Petersen.
\newblock {Directional Distance Functions in DEA with Optimal Endogenous
  Directions}.
\newblock \emph{Operations Research}, 66\penalty0 (4):\penalty0 1068--1085,
  2018.

\bibitem[Rigby and Bilodeau(2015)]{ManagementToolsTrends}
D.~Rigby and B.~Bilodeau.
\newblock Management tools and trends 2015.
\newblock Bain \& Company, 2015.

\bibitem[Rigby(2015)]{ManagemenTools2015}
D.K. Rigby.
\newblock Management tools 2015 - an executive's guide.
\newblock Bain \& Company, 2015.

\bibitem[Rostami et~al.(2017)Rostami, Neri, and Epitropakis]{Rostamy_ea2017}
S.~Rostami, F.~Neri, and M.~Epitropakis.
\newblock Progressive preference articulation for decision making in
  multi-objective optimisation problems.
\newblock \emph{Integrated Computer-Aided Engineering}, 24\penalty0
  (4):\penalty0 315--335, 2017.

\bibitem[Rudin et~al.(2022)Rudin, Chen, Chen, Huang, Semenova, and
  Zhong]{rudin2022interpretable}
C.~Rudin, C.~Chen, Z.~Chen, H.~Huang, L.~Semenova, and C.~Zhong.
\newblock Interpretable machine learning: Fundamental principles and 10 grand
  challenges.
\newblock \emph{Statistics Surveys}, 16:\penalty0 1--85, 2022.

\bibitem[Schaffnit et~al.(1997)Schaffnit, Rosen, and Paradi]{SchaffnitEA1997}
C.~Schaffnit, D.~Rosen, and J.C. Paradi.
\newblock Best practice analysis of bank branches: An application of {DEA} in a
  large {Canadian} bank.
\newblock \emph{European Journal of Operational Research}, 98\penalty0
  (2):\penalty0 269--289, 1997.

\bibitem[Silva~Portela et~al.(2003)Silva~Portela, Borges, and
  Thanassoulis]{Portela2003}
M.C. Silva~Portela, P.C. Borges, and E.~Thanassoulis.
\newblock Finding closest targets in non-oriented {DEA} models: The case of
  convex and non-convex technologies.
\newblock \emph{Journal of Productivity Analysis}, 19:\penalty0 251--269, 2003.

\bibitem[Thach(1988)]{thach1988design}
P.T. Thach.
\newblock The design centering problem as a {DC} programming problem.
\newblock \emph{Mathematical Programming}, 41\penalty0 (1):\penalty0 229--248,
  1988.

\bibitem[Wachter et~al.(2017)Wachter, Mittelstadt, and
  Russell]{wachter2017counterfactual}
S.~Wachter, B.~Mittelstadt, and C.~Russell.
\newblock Counterfactual explanations without opening the black box: Automated
  decisions and the {GDPR}.
\newblock \emph{Harvard Journal of Law \& Technology}, 31:\penalty0 841--887,
  2017.

\bibitem[Zhu(2016)]{zhu_book2016}
J.~Zhu.
\newblock \emph{Data Envelopment Analysis - A Handbook on Models and Methods}.
\newblock Springer New York, 2016.

\bibitem[Zofio et~al.(2013)Zofio, Pastor, and
  Aparicio]{ZofioPastorAparicio2013}
J.~L. Zofio, J.~T. Pastor, and J.~Aparicio.
\newblock The directional profit efficiency measure: on why profit inefficiency
  is either technical or allocative.
\newblock \emph{Journal of Productivity Analysis}, 40:\penalty0 257--266, 2013.

\end{thebibliography}

\section*{Appendix}

\renewcommand{\thesubsection}{\Alph{subsection}}

\setcounter{equation}{29}
\renewcommand\theequation{\arabic{equation}}

\label{sec:extensions}

Here, we extend the analysis in Section \ref{sec:bilevel} by investigating alternative returns to scale and by investigating changes in the outputs rather than the inputs. In both extensions, we will consider a combination of the $\ell_0$, $\ell_1$ and $\ell_2$ norms as in objective function \eqref{eq:costf}.

\subsection{Changing the returns to scale}

In Section \ref{sec:bilevel}, we have considered the DEA model with constant return to scale (CRS), where the only requirement on the values of $\bm{\lambda}$ is that they are positive, i.e., $\bm{\lambda} \in \mathbb{R}_+^{K+1}$, but we could consider other technologies. In that case, to be able to transform our bilevel optimization problem to a single-level one, we should take into account that for each new constraint derived from the conditions on $\bm{\lambda}$, a new dual variable has to be introduced. We will consider the varying return to scale (VRS) model as it is one of the most preferred one by firms \citep{Bogetoft2012}, but extensions to other models are analogous.

Consider the input case. With the same transformation as before, we have:
\begin{align*}
	\footnotesize
	\min_{\bm{\hat x},F} \quad & C(\bm{x}^0,\bm{\hat x}) \\
	\text{s.t.} \quad & \bm{\hat x} \in \mathbb{R}^I_{+} \\
	&F \leq F^{*} \\
	& F \in \argmin_{\bar F, \lambda^0, \dots, \lambda^K} \footnotesize \left\{ \right. \bar F :
	\quad   {\bm{\hat x}} \geq \sum_{k=0}^K \beta^k \bm{x}^k,
	\hat F\bm{y}^0 \leq \sum_{k=0}^K \beta^k \bm{y}^k, \bar F \geq 0,  \bm{\beta} \in \mathbb{R}^{K+1}_{+}, \sum_{k=0}^K \beta^k = F \left. \footnotesize \right\}.
\end{align*}

Notice that the only difference is that we have a new constraint associated with the technology, namely, $\sum_{k=0}^K \beta^k = F$. Let $\kappa \geq 0$ be the new dual variable associated with this constraint. Then, the following changes are made in constraints \eqref{eq:stattf} and \eqref{eq:stattf2}:
\begin{align}
	\bm{\gamma}_{O}^\top\bm{y}^0&+\kappa = 1 \label{eq:stattvrs}\\
	\bm{\gamma}_I^\top \bm{x}^k - \bm{\gamma}_{O}^\top \bm{y}^k& -\kappa\geq  0 \quad k=0,\dots,K \label{eq:statt2vrs}.
\end{align}

The single-level formulation for the counterfactual problem for VRS DEA is as follows:
\begingroup
\allowdisplaybreaks
\begin{align}
	\tag{CEVDEA}
	\min_{\bm{\hat x},F, \bm{\beta},\bm{\gamma}_I,\bm{\gamma}_O,\bm{u},\bm{v},\bm{w},\kappa ,\bm{\eta},\bm{\xi}} \quad & \nu_0 \sum_{i=1}^I \xi_i +\nu_1 \sum_{i=1}^I \eta_i + \nu_2  \sum_{i=1}^I (x_i^0-\hat{x}_i)^2 \nonumber \\
	\text{s.t.} \quad 	&F\leq F^{*}\nonumber\\
	&\sum_{k=0}^K \beta^k = F\nonumber \\
	&\bm{\hat x} \in \mathbb{R}^I_{+}\nonumber \\
	&	\kappa\geq 0 \nonumber\\
	&	\bm{u},\bm{v},\bm{w}  \in \{0,1\} \nonumber\\
	&\eqref{eq:primm}-\eqref{eq:primm_l},\eqref{eq:M1}-\eqref{eq:M4}\nonumber\\
	&\eqref{eq:stattf3},\eqref{eq:l0}-\eqref{eq:statt2vrs}.\nonumber
\end{align}
\endgroup

\subsection{Changing the outputs}

We have calculated the counterfactual instance of a firm as the mininum cost changes in the inputs in order to have a better efficiency. In the same vein, we could consider instead changes in the outputs, leaving the same inputs. Again, suppose firm 0 is not fully efficient, $E^{0} <1$. Now, we are interested in calculating the minimum changes in the outputs $ \bm y^0$ that make it to have a higher efficiency $E^{*}>E^{0}$. Let $\hat {\bm y}$ be the new outputs of firm 0 that make it to be at least $E^{*}$ efficient. We have, then, the following bilevel optimization problem:
\begin{align*}
	\footnotesize
	\min_{\bm{\hat y},E} \quad & C(\bm{y}^0,\bm{\hat y}) \\
	\text{s.t.} \quad & \bm{\hat y} \in \mathbb{R}^O_{+} \\
	&E \geq E^{*} \\
	& E \in \argmin_{\bar E, \lambda^0, \dots, \lambda^K} \footnotesize \left\{ \right. \bar E :
	\quad\bar  E \bm{x}^0\geq \sum_{k=0}^K \lambda^k \bm{x}^k,
	\hat{\bm{y}} \leq \sum_{k=0}^K \lambda^k \bm{y}^k, \bar E \geq 0,  \bm{\lambda} \in \mathbb{R}^{K+1}_{+} \left. \footnotesize \right\}.
\end{align*}

Following similar steps as in previous section, the single-level formulation for the counterfactual problem in DEA in the output case is as follows:
\begingroup
\allowdisplaybreaks
\begin{align}
	\tag{CEODEA}
	\min_{\bm{\hat y},E, \bm{\lambda},\bm{\gamma}_I,\bm{\gamma}_O, \bm{u},\bm{v},\bm{w},\bm{\eta},\bm{\xi}} \quad & \nu_0 \sum_{o=1}^O \xi_o +\nu_1 \sum_{o=1}^O \eta_o + \nu_2  \sum_{o=1}^O (y_o^0-\hat{y}_o)^2\\ \nonumber\\
	\text{s.t.} \quad & \bm{\hat y} \in \mathbb{R}^O_{+} \nonumber\\
	&E\geq E^{*}\nonumber\\
	&E \bm{x}^0\geq \sum_{k=0}^K \lambda^k \bm{x}^k \nonumber \\
	& \bm{\hat y} \leq \sum_{k=0}^K \lambda^k \bm{y}^k \nonumber\\
	&\bm{\gamma}_{I}^\top \bm{x}^0= 1\nonumber \\
	&\bm{\gamma}_{O}^\top \bm{y}^k - \bm{\gamma}_{I}^\top \bm{x}^k \leq 0 && k=0,\dots, K\nonumber\\
	&\gamma^i_I \leq M_I u_i \quad&& i=1,\dots,I\nonumber\\
	& E x_i^0-\sum_{k=0}^K \lambda^k x_i^k \leq M_I(1-u_i) \quad &&i=1,\dots,I \nonumber\\
	& \gamma^o_{O} \leq M_O v_o \quad&& o=1,\dots,O\nonumber\\
	& -\hat{y}_o+\sum_{k=0}^K \lambda^k y_o^k \leq M_O(1-v_o) \quad &&o=1,\dots,O\nonumber\\
	&\lambda^k\leq M_f w_k \quad&& k=0,\dots,K\nonumber\\
	& \bm{\gamma}_{O}^\top \bm{y}^k - \bm{\gamma}_{I}^\top \bm{x}^k \leq M_f(1-w_k) \quad&& k=0,\dots, K\nonumber\\
	&	-M_{\text{zero}}\xi_o \leq y_o^0-\hat{y}_o && o=1,\dots,O\nonumber\\
	& y_o^0-\hat{y}_o \leq M_{\text{zero}} \xi_o && o=1,\dots,O\nonumber \\
	&	\eta_o \geq y_o^0-\hat{y}_o&& o=1,\dots,O\nonumber\\
	&-\eta_o \leq y_o^0-\hat{y}_o && o=1,\dots,O \nonumber\\
	& E,\bm{\lambda},\bm{\gamma}_I,\bm{\gamma}_O,\bm{\eta} \geq 0 \nonumber\\
	&\bm{u},\bm{v},\bm{w},\bm{\xi}  \in \{0,1\}.\nonumber
\end{align}
\endgroup
As in the input model, depending on the cost function, we either obtain an MILP model or a Mixed Integer Convex Quadratic model with linear constraints. This model could be formulated analogously for the VRS case.

\end{document}